\documentclass[a4paper,12pt]{article}

\usepackage[latin1]{inputenc}
\usepackage[T1]{fontenc}
\usepackage{latexsym,amssymb,amsmath,epsfig}
\usepackage{a4wide}
\usepackage{graphicx}

\newcommand \eps {\varepsilon}
\renewcommand \leq {\leqslant}
\renewcommand \geq {\geqslant}

\newtheorem{defi}{Definition}[section]
\newtheorem{lemm}[defi]{Lemma}

\newtheorem{theo}[defi]{Theorem}
\newenvironment{dem}{\vskip 2mm\noindent {\it Proof}:}
                    {\hfill $\square$ \vskip 2mm \noindent} 
\newenvironment{demof}[1]{\vskip 2mm\noindent {\it Proof of #1}:}
                    {\hfill $\square$ \vskip 2mm \noindent} 
 
\newenvironment{demoof}[1]{{ \em \bf {\it  Proof of #1: }}}{$\hfill \Box$ \smallskip}

\title{Image denoising by statistical area thresholding}
\author{D. Coupier, A. Desolneux, and B. Ycart} 

\begin{document}
\maketitle
\begin{center}
MAP5, UMR CNRS 8145, Universit\'e Ren\'e Descartes, Paris
\end{center}
\vskip 0.3cm
{\bf E-Mail addresses :} {\tt
  \{coupier,desolneux,ycart\}@math-info.univ-paris5.fr}\\  

\vskip 0.3cm 
\noindent
{\bf Corresponding author :} A. Desolneux

\noindent
{\bf Mail address :} MAP5, UFR Math-Info,\\
\hspace*{2.8cm} Universit\'e Paris 5,  45 rue des Saints-P\`eres\\
\hspace*{2.8cm} 75270 PARIS Cedex 06, FRANCE.

\noindent
{\bf Telephone :} 33 1 44 55 35 26 \\
{\bf Fax :} 33 1 44 55 35 35
\vskip 2cm

\begin{abstract}
Area openings and closings are morphological filters which efficiently suppress
impulse noise from an image, by removing small connected components of level
sets. The problem of an objective choice of threshold for the area
remains open. Here, a mathematical model for random images will be
considered. Under this model, a Poisson approximation for
the probability of appearance of any local
pattern can be computed. In particular, the probability of observing a 
component with size larger than $k$ in pure impulse noise has an explicit
form. This permits the definition of a statistical test on the significance of
connected components, thus providing an explicit formula for the area 
threshold of the denoising filter, as a function of the impulse noise
probability parameter. Finally, using threshold decomposition, a
denoising algorithm for grey level images is proposed. 

\end{abstract}
\vskip 1cm


\newpage

\section{Introduction}
\label{intro}

The general problem of image denoising consists of deciding what is the
``signal'' and should be kept, and what is the noise, and
must be removed. Many different criteria can be used to detect 
the noise-induced structures. For example, the oscillations due to an additive
gaussian noise can be measured in terms of the  
wavelet coefficients. The noise may then be removed by a thresholding 
in the wavelet domain. Donoho and Johnstone \cite{donoho} gave an 
explicit way to choose the threshold as a function of the variance 
of the noise. Their claim is that ``denoising, with high probability, 
rejects pure noise completely''. The underlying idea is that in 
pure noise, all the structures that actually belong to the image
could not appear; or else, the structures coming from the image itself
can be defined as those ``objects'' which would have a very small 
probability of appearing in a pure noise. This idea was implemented 
in \cite{Desolneux} and \cite{DMM2}  for the detection of
alignments and meaningful level lines in an image.

Here, we shall focus on the size of connected components
of the level sets of the image. Removing small components is a classical
and efficient way of removing impulse noise from an image. This method, 
known as ``the grain filter'', was first introduced  in the framework of 
Mathematical Morphology \cite{Serra} by Vincent in \cite{vincent} 
as morphological area openings and closings (see also \cite{vincentbis} and
\cite{Heijmans}). This filter is sometimes  
called the ``extrema killer''. It was then generalized by Masnou and 
Morel in \cite{masnou98}, and by Monasse and Guichard in \cite{monasse}. In
\cite{sapoval95} a similar filter was used, in the framework
of gradient percolation, for recovering fuzzy images.

But the main question remains: how should the threshold for the area of the
components that have to be kept, be chosen? A natural idea, imported 
from statistical inference, consists of fixing an a priori risk
level $\eps$ (e.g. $\eps = 0.001$), and deciding that anything that 
has probability lower than $\eps$ of occurring under a pure noise 
hypothesis cannot come from the noise and hence should be kept in the
image. Thus for the threshold area, one will choose the integer 
$s(n,p,\eps)$,  such that a connected component of size 
$k\geq s(n,p,\eps)$ has a probability less than $\eps$ of appearing in a 
pure noise image with probability parameter $p$ and size $n\times n$.
Applying a grain filter with area threshold $s(n,p,\eps)$ will ensure
that, with probability larger than $1-\eps$, pure noise is eliminated.
To implement this, one must be able to compute the probability of
any connected component of size $k$ appearing in a pure noise
image. An exact computation is not feasible. However
an approximation can be given if the image is large: our main
theoretical result (Theorem \ref{th:psi}) gives a Poisson
approximation for the probability of occurrence for 
any image property which is local in the sense that its definition 
involves only a fixed number of connected pixels.
\vskip 3mm
Our plan is as follows. Section \ref{part_pa.sec} is devoted to
the probabilistic model of noise in binary images: all pixels are
independent, black with probability $p$ or white with probability
$1\!-\!p$. The Poisson approximation result  will be stated (Theorem
\ref{th:psi}) and an
outline of its proof will be given (technical details will be postponed
to the Appendix). Section \ref{part_denois.sec} is devoted to 
applications. We will first explain
how Theorem \ref{th:psi}, together with numerical combinatorial results on square 
lattice animals\footnote{Square lattice animals or ``polyominoes'' are simply
  defined as 
  connected clusters of squares in the plane (for example, the ``Tetris'' game
  uses all lattice animals of size $4$).}, can be used to obtain an explicit
formula for the  
size threshold $s(n,p,\eps)$. An example of denoising for a binary
image will be given. Then we shall extend the method to grey level
images through threshold decomposition: the binary image corresponding
to each grey level is treated separately, then all denoised binary images are
recombined. Some experiments and a discussion of the obtained results
come last.

\section{Probability of a local property}
\label{part_pa.sec}

Our probabilistic model for random images is the following.
Let $n$ be a positive integer. 
Consider the \textit{pixel set} $\Xi_{n} = \lbrace 1, 
\ldots ,n \rbrace ^{2}$. 
A \textit{binary image of size} $n$ is a mapping from $\Xi_n$ to 
$\lbrace 0,1 \rbrace$ (black/white). 
Their set is denoted by $E_{n}$. It is endowed with 
the probability distribution 
$\mu_{n,p}$ defined by: 
each pixel is black with probability $p$ or white with probability
$1-p$, and all the pixel colors are independent.
A \textit{random image of size} $n$ \textit{and probability parameter} 
$p$, denoted by $\mathcal{I}_{n,p}$, is a random element of $E_{n}$ 
with probability distribution 
$\mu_{n,p}$.

\begin{figure}[h]
\begin{center}
\epsfxsize=5cm 
\epsfbox{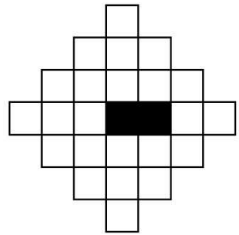} 
\end{center}
\caption{Example of an image on the ball $B(0,r)$ with $r=3$. This small image
  is also called a {\it pattern}. The number of black pixels of this pattern
  $D$ is $b(D)=2$.}
\label{ball.fig}
\end{figure}

The pixel set $\Xi_{n}$ is embedded in $\mathbb{Z}^{2}$ and 
naturally endowed with a \textit{graph structure}. 
We consider in this paper the case of $4$-connectivity ($2$ horizontal 
and $2$ vertical neighbors). 
For purely technical reasons, it will be convenient that all pixels
have the same neighborhood: this is why we impose periodic boundary
conditions, deciding that $(1,j)$ is a neighbor of $(n,j)$ and $(j,1)$ 
of $(j,n)$. Thus the graph is a $2$-dimensional torus. 
As usual, the \textit{graph distance} $d$ is defined as the minimal 
length of a path between two pixels. 
We shall denote by $B(x,r)$ the ball of center 
$x$ and radius $r$ with respect to the distance $d$. It is defined by
$$
B(x,r) = \lbrace y \in \Xi_{n}; \; d(x,y) \leq r \rbrace ~.
$$
Notice that this ball $B(x,r)$ is diamond-shaped (it is a rhombus) and that for
$r<n/2$, it contains $2r^{2}+2r+1$ pixels (see Figure \ref{ball.fig}). 
For the rest of this section, the radius $r$ is a fixed integer, and
the image size $n$ is larger than $2r+1$.

The image properties we are interested in are all local, in the sense
that they can be described inside balls of a fixed radius.
All balls are translations of each other. We shall choose a
ball of radius $r$, say $B(0,r)$, and fix a translation
$\tau_x$, from $B(0,r)$ to $B(x,r)$ for all $x$.
We call {\it pattern}, and denote by $D$, an image defined on $B(0,r)$, 
and determined by its set of black pixels, denoted by
$\beta(D)$ (see Figure \ref{ball.fig} for an example of pattern). Of course,
  $B(0,r)\setminus\beta(D)$ is the set of  
white pixels. We shall denote by $b(D)$ the cardinality of $\beta(D)$ 
(number of black pixels in the pattern).  We shall deal with
rather small levels of noise, seen as relatively sparse black 
pixels on a white background. This is of course a mere convention: 
swapping black and white, together with $p$ and $1\!-\!p$ 
does not change the model. Thus, in what follows, we will always assume that
$p\leq \frac{1}{2}$.

If $D$ is a pattern on $B(0,r)$ and $\tau$ is a translation of pixels, 
we shall denote by $\tau(D)$ the pattern on $B(\tau(0),r)$, whose set
of black pixels is $\tau(\beta(D))$. If $\tau(0)=x$, we denote by $D(x)$ 
the property: ``the restriction of the image to $B(x,r)$ is 
$\tau(D)$''. The property we are actually interested in is
$$
\widetilde{D} = (\exists x\in \Xi_{n} \,,\; D(x))\;.
$$
In other words $\widetilde{D}$ means: ``a copy of pattern $D$ can be found
somewhere in the image''.
\vskip 3mm
The patterns $D$ are the building blocks of all local
properties. Indeed, there exists only a finite number of such patterns
(precisely $2^{2r^2+2r+1}$): let us denote their set by $\cal{D}$.
Any {\it assertion} relative to the pixels in $B(0,r)$ will be called
``local'':  it can be expressed in a unique way as a disjunction
(logical ``or'', denoted by $\vee$) of distinct patterns. 

The following definitions will be used in the counting 
of occurrences of a local property in an image. 
\begin{defi}
\label{de:localassertion}
Let $\psi$ be a local assertion, relative to the pixels in $B(0,r)$. 
\begin{enumerate}
\item
The {\rm definition set} of $\psi$, denoted by ${\cal D}(\psi)$, is
the subset of ${\cal D}$ such that
$$
\psi = \bigvee_{D \in \mathcal{D}(\psi)} D ~.
$$
\item
The {\rm black index} $b(\psi)$ of $\psi$ is the integer 
$b(\psi)$ defined by 
$$
b(\psi) = \min_{D \in {\cal D}(\psi)} \{ b(D) \} ~.
$$
\item
A {\rm meaningful definition set} of $\psi$,  
denoted by $\mathcal{D}_0(\psi)$, is a subset of ${\cal D}(\psi)$
such that
\begin{enumerate}
\item
$
\forall D\in {\cal D}_0(\psi)\;,\quad b(D)=b(\psi)\;,
$
\item
If $\tau$ is a translation, then
$
D,D'\in {\cal D}_0(\psi)\mbox{ and } \tau(\beta(D))=\beta(D') \mbox{ imply } D=D'\;,
$
\item
$
D\in {\cal D}(\psi) \mbox{ and } b(D)=b(\psi) \mbox{ imply }
\exists \tau\,, \exists D'\in {\cal D}_0(\psi) \, , \mbox{ s. t. }
\tau(\beta(D))=\beta(D')\;.  
$
\end{enumerate}
All meaningful definition sets have the same cardinality, which will
be called the {\rm meaningful index} of $\psi$, and denoted by $e(\psi)$.
\end{enumerate}
\end{defi}
The black index $b(\psi)$ is the minimal number of black pixels, in a
pattern that satisfies $\psi$. One can see the meaningful index 
$e(\psi)$ as the maximal number of patterns with exactly $b(\psi)$ 
black pixels that satisfy $\psi$, up to possible translations. 
Both will be used to count occurrences of the {\it local property} 
based on $\psi$.   

\smallskip

\noindent {\bf Example.} Let us illustrate all these definitions by
considering a simple example: the property ``there exist two
connected black 
pixels''. On the ball of radius $r=1$, the definition set of this local
assertion $\psi$ is composed
of all those patterns on $B(0,1)$ whose center is black, and
at least one of the $4$ neighbors is also black ($15$ patterns). The black
index $b(\psi)$ of $\psi$ equals $2$, and its meaningful index is $e(\psi)=2$
(a possible meaningful definition set is made of the two
 patterns on  $B(0,1)$ such that the center and its right horizontal
neighbor, resp. its top vertical neighbor, are the only black pixels in the
ball $B(0,1)$).  

\smallskip

\begin{defi}
Let $\psi$ be a local assertion, and $\psi(x)$ its localization 
on the ball centered at $x$~:
$$
\psi(x) = \bigvee_{D \in \mathcal{D}(\psi)} D(x) ~.
$$
We call {\rm local property} based on $\psi$, and denote by 
$\widetilde{\psi}$ the property 
$$
\widetilde{\psi} = (\exists x\,,\; \psi(x))\;.
$$ 
\end{defi}
Our basic example of a local property $\widetilde{\psi}$ 
is: ``there exists a connected
component of $k$ black pixels''. A connected component of size $k$ 
is always included in a ball of radius $r\geq k/2$. 
The local assertion $\psi$ is 
``there exists a connected component of
size $k$ in $B(0,r)$''. The definition set is
the set of all patterns on $B(0,r)$ having at least $k$ connected black 
pixels. The black index is the minimal number of black pixels necessary for
$\psi$ to be satisfied (obviously $k$ in our example). The meaningful
index is the number of connected components of size $k$, up to
translations (see Section \ref{part_denois.sec}). 
\vskip 3mm
For a fixed level $p$ with $0<p<1$, if we let $n$ tend to infinity, by the 
independence of pixels, it is easy to see that 
asymptotically any pattern will be 
present in a random image with a probability tending to $1$ 
(see \cite{CDY03} for more precise results). Therefore the asymptotic probability
for the random image ${\cal I}_{n,p}$ to satisfy 
$\widetilde{\psi}$ is $1$, whatever $\psi$. 
That asymptotic probability can be different from $1$ only if
$p=p(n)$ tends to $0$ as $n$ tends to infinity.
Thus our images will have a relatively small proportion of black
pixels.
\vskip 3mm
A classical object of the theory of random graphs (see
\cite{Bollobas,Spencer} as general references),  is the notion 
of \textit{threshold function}. It describes 
the appearance of a given subgraph in a random graph. The notion of
threshold function easily adapts to random images.
Let $\mathcal{A}$ be an image property. 
The function $\theta(n)$ is called a threshold function of $\mathcal{A}$ 
if for $p(n)\leq 1/2$ then
$$
\lim_{n \to \infty} \frac{p(n)}{\theta(n)} = 0 
\Longrightarrow
\lim_{n \to \infty} \mu_{n,p(n)} (\mathcal{A}) = 0\;,
$$ 
and 
$$
\lim_{n \to \infty} \frac{p(n)}{\theta(n)} = \infty
\Longrightarrow 
\lim_{n \to \infty} \mu_{n,p(n)} (\mathcal{A}) = 1\;.
$$ 
Notice that a threshold function is not unique.  
For instance if $\theta(n)$ is a threshold function for $\mathcal{A}$, 
then so is $c\theta(n)$ for any positive constant $c$. It is customary 
to ignore this and talk about ``the'' threshold function of
$\mathcal{A}$.
We then have the following lemma.
\begin{lemm}
\label{le:threshold}
The threshold function of the local property $\widetilde{\psi}$ is 
$n^{-\frac{2}{b(\psi)}}$.
\end{lemm} 
\begin{dem}
We shall just give here the main steps, since the detailed proof will appear in
\cite{CDY03}. Let $D$ be a pattern and let $X_{n}$ denote the number of copies
of $D$ in the image. Then 
$$\mu_{n,p(n)}(\tilde{D})= \mu_{n,p(n)}(X_{n}>0) \leq \mathbb{E}(X_{n}) =
n^{2}p(n)^{b(D)}(1-p(n))^{2r^2+2r+1-b(D)} \leq
n^{2}p(n)^{b(D)} .$$
On the other hand, let $Y_n$ denote the number
of copies of $D$ occurring in balls $B(x,r)$ where both coordinates of $x$ are
multiples of $2r+1$ (which implies that two such balls cannot meet). The number
of such balls is $n_r^2$ where $n_r=\lfloor \frac{n}{2r+1} \rfloor$.
Then
\begin{eqnarray*}
\mu_{n,p(n)}(X_{n}>0) \geq  \mu_{n,p(n)}(Y_{n}>0) & = &
1-\mu_{n,p(n)}(Y_{n}=0) \\
& = & 1-\left( 1-p(n)^{b(D)}(1-p(n))^{2r^2+2r+1-b(D)}\right)^{n_r^2} \\
& \geq & 1- \exp(-{n_r}^2 p(n)^{b(D)}(1-p(n))^{2r^2+2r+1-b(D)}) 
\end{eqnarray*}
Using these inequalities and the definition of a threshold function, we conclude that 
$\theta(n)=n^{-2/b(D)}$ is the threshold function of the property
$\tilde{D}$ . To conclude, one has to check that the
threshold function of a disjunction of patterns is
the smallest threshold function of these
patterns.
\end{dem}

Lemma \ref{le:threshold} means that the appearance of a local property
mainly depends on its 
black index: if $p(n)$ is small compared to 
$n^{-\frac{2}{b}}$, then the probability of any 
local property that needs $b$ black pixels to be satisfied is
small. If $p(n)$ is large compared to $n^{-\frac{2}{b}}$, then the
probability is large.
The particular case $b(\psi) = 0$ corresponds to the appearance of a white
ball. If there exists $\alpha$, with $0\leq\alpha<1$, such that for all $n$,
we have $p(n) \leq \alpha$, then the probability for a white ball of being
present in 
the random image always tends to $1$ as $n$ tends to infinity: there is no
threshold function. 
From now on, we will always assume that the black index of 
$\psi$ is positive.
\vskip 3mm
Lemma \ref{le:threshold} suggests that the correct scaling for
$p(n)$ when one studies a local property $\widetilde{\psi}$ is 
$p(n) = c\,n^{-\frac{2}{b(\psi)}}$. Our main result shows that with
this scaling, the probability of $\widetilde{\psi}$ in a random image
converges to a non trivial limit.
\begin{theo}
\label{th:psi}
Let $\psi$ be an assertion on $B(0,r)$, with black index $b(\psi)$ and
meaningful index $e(\psi)$.  
Let 
$p(n) = cn^{-\frac{2}{b(\psi)}}$, where $c$ is a positive constant.
Then
\begin{equation}
\label{limth}
\lim_{n \to \infty} \mu_{n,p(n)} (\widetilde{\psi}) 
= 1-\exp( - e(\psi) c^{b(\psi)}) ~.
\end{equation}
\end{theo}
The reason why such a result is called a Poisson approximation becomes
clear if one considers the property ``there exists a black
pixel''. Let $X_n$ be the total number of black pixels.
Since all pixels are independent, the random variable $X_n$ follows the
binomial distribution with parameters $n^2$ and $p(n)$. In particular
the probability that there exists a black pixel is
$$
\mathbb{P}[X_n>0] = 1-(1-p(n))^{n^2}\;.
$$
Here the black index is $1$ and the threshold function is $n^{-2}$. Take
$p(n)=cn^{-2}$. Then the binomial distribution of $X_n$ converges to the
Poisson distribution with parameter $c$, and the probability that there
exists a black pixel ($X_n>0$) tends to $1-\exp(-c)$. 

The situation is not so
simple as soon as the black index is larger than $1$. Consider for
instance again the local property $\widetilde{\psi}$: ``there exist two
connected black 
pixels''. We already saw that on the ball of radius $r=1$, the definition set
is composed 
of all those patterns on $B(0,1)$ whose center is black, and
at least one of the $4$ neighbors is also black ($15$ patterns).
Consider the number of occurrences of any of those
patterns, somewhere in the random image. It is a sum of
Bernoulli random variables. However
they are not independent: patterns on balls centered
at two adjacent pixels have one pixel in common. The same can be said
of any local property $\widetilde{\psi}$: the number of occurrences of
$\psi(x)$ can be viewed as a sum of (dependent) Bernoulli random
variables. The sum of a large number of Bernoulli r.v.'s converges in
distribution to a Poisson distribution, provided the dependencies
between the variables are not too large. In the theory of random graphs, 
similar results are frequent (see e.g. \cite{Spencer} Lecture 1 p.296, Lecture
2 p.303 or Lecture 5 p.314).  

\begin{demof}{Theorem \ref{th:psi}}
There are several ways to prove a Poisson approximation result. 
We chose the famous ``moment method'' based on the following result
(\cite{Bollobas}, Chapter 1 p.25). 
\begin{lemm}
\label{moments}
Let $(X_{n})_{n \in \mathbb{N}^{\ast}}$ be a sequence of integer valued, 
nonnegative random variables and 
$\lambda$ be a strictly positive real. 
For all $n , l \in \mathbb{N}^{\ast}$ define the quantity
$$
E_{l}(X_{n}) = \sum_{k \geq l} \mathbb{P}(X_{n} = k)\frac{k!}{(k-l)!}
~.
$$
If, for all $l \in \mathbb{N}^{\ast}$, 
$\lim_{n \to \infty} E_{l}(X_{n}) = \lambda^{l}$ then 
$(X_{n})$ converges in distribution to 
the Poisson distribution with parameter $\lambda$. 
\end{lemm}
In our case, $X_n$ counts the number of occurrences in the
random image of some patterns, to be precised later. The ``moment''
$E_l(X_n)$ is the expected number of 
ordered $l$-tuples of occurrences of those patterns.
\vskip 3mm
Firstly, one should observe that
patterns in the definition set of $\psi$ cannot be all treated
equally: since $p(n)=cn^{-\frac{2}{b(\psi)}}$, by Lemma
\ref{le:threshold} any pattern with more
than $b(\psi)$ black pixels has a vanishing probability of being
observed. Hence we can reduce the set of patterns to those having
exactly $b(\psi)$ black pixels. In the example of two connected
pixels with $r=1$, ${\cal D}(\psi)$ has $15$
different patterns, but only $4$ of them have exactly $2$ black pixels.

Now, one has to take care of multiple counts. Among the $4$ patterns
on $B(0,1)$ that have $2$ black pixels, $2$ patterns have two
horizontal black neighbors, and the $2$ other patterns have two
vertical black neighbors. Assume the image has only one occurrence
of two horizontal black neighbors. If we examine all possible pixels 
$x$, we will find two adjacent centers for which $\psi(x)$ is
satisfied. In order
to obviate this problem, we need to count patterns {\it up to possible
translations}. We say that two patterns with black index $b(\psi)$ are
equivalent if their sets of black pixels are translations of each
other. The number of equivalence classes is the meaningful index 
$e(\psi)$ of Definition \ref{de:localassertion}. 
(In the example of two connected black pixels,
there are two equivalence classes: horizontal or vertical neighbors).

We choose a meaningful set, i.e. we 
fix a pattern for each equivalence class: 
$$
{\cal D}_0(\psi) = \{\bar{D}_{1}, \ldots , \bar{D}_{e(\psi)}\}~.
$$
The counting variable $X_n$ to which Lemma \ref{moments} will be
applied is the total number of occurrences of one of the patterns 
$\bar{D}_{1}, \ldots , \bar{D}_{e(\psi)}$, in the random image 
${\cal I}_{n,p(n)}$:
$$
X_n = \sum_{x\in \Xi_n} \sum_{i=1}^{e(\psi)}
\mathbb{I}_{\bar{D}_i(x)}({\cal I}_{n,p(n)})\;,
$$
where $\mathbb{I}$ denotes the indicator function of an event. 
The expectation of $X_n$ is
$$
\mathbb{E}(X_n) = n^2\,e(\psi)\, (p(n))^{b(\psi)}(1-p(n))^{2r^2+2r+1-b(\psi)}
~.
$$
As $n$ tends to infinity, it tends to $e(\psi)\,c^{b(\psi)}$, 
which is the parameter of the Poisson approximation in formula 
(\ref{limth}). In order to apply Lemma \ref{moments} to $X_n$, one
has to check that the hypothesis holds. 
\begin{lemm}
\label{hypmoments}
$$
\forall l \in \mathbb{N}^{\ast} \, , \hspace{0.5cm} \lim_{n\to\infty} E_l(X_n) = \left( e(\psi) c^{b(\psi)} \right)^l\;.
$$
\end{lemm}
The proof of Lemma \ref{hypmoments} is rather
technical and will be given in the Appendix.
\vskip 3mm
Now Lemma \ref{moments} implies that $X_n$ converges in
distribution to the Poisson distribution with parameter 
$e(\psi) c^{b(\psi)}$. Therefore $\mu_{n,p(n)}(X_n>0)$ tends to
$1-\exp(-e(\psi) c^{b(\psi)})$.
It is clear that $X_n>0$ implies that ${\cal I}_{n,p(n)}$ satisfies
$\widetilde{\psi}$. Hence 
$\mu_{n,p(n)}(X_n>0)\leq \mu_{n,p(n)}(\widetilde{\psi})$. 
Moreover, the event $(\widetilde{\psi} \setminus (X_n>0))= (\widetilde{\psi} \cap
(X_n=0))$ implies the appearance of a pattern with at least $b(\psi)+1$ black
pixels in a ball of radius $r$, and by Lemma \ref{le:threshold}, its
probability tends to $0$ as $n$ tends to infinity.  
Therefore, 
$$
\lim_{n\to\infty} \mu_{n,p(n)}(X_n>0) = 
\lim_{n\to\infty} \mu_{n,p(n)}(\widetilde{\psi}) =
1-\exp( - e(\psi) c^{b(\psi)})
\;.
$$ 
It should be noticed that the asymptotics of $X_n$ does not depend on
the choice of the meaningful definition set $\{\bar{D}_{1}, \ldots ,
\bar{D}_{e(\psi)}\}$. It does not depend either on the radius $r$ of the
ball. Consider for instance the property $\widetilde{\psi}$
``the image contains two horizontally connected black pixels''. 
Its definition set for the ball $B(0,r)$ has $r^2 2^{2r^2+2r}$
elements. Among these, only $2r^2$ have exactly $2$ black pixels, and there is only
one equivalence class up to translations, whatever $r$. Therefore $r$
is a phantom parameter, as should be expected. It serves only to
ensure that properties remain local.
\end{demof}

\section{Application to image denoising}
\label{part_denois.sec}

In the previous section, we computed the asymptotic probability of
appearance of any local property in a random binary image. This provides
the basis of a statistical test to decide whether an observed pattern in an image
may be due to 
noise or not, and this test can be applied for image denoising.
In this section, all the considered images will be corrupted by 
the same kind of noise, namely impulse noise. 
This type of noise models for example the fact that some (unknown) part of the
data is lost. We will assume that the probability parameter of the noise is known. 
We will first start with the denoising of binary images, and then extend it to
grey level images using their threshold decomposition.

\subsection{Binary images}
\label{part_bin.subsec}

\begin{figure}[h]
\begin{center}
\begin{tabular}{cccc}
\epsfxsize=4cm 
\epsfbox{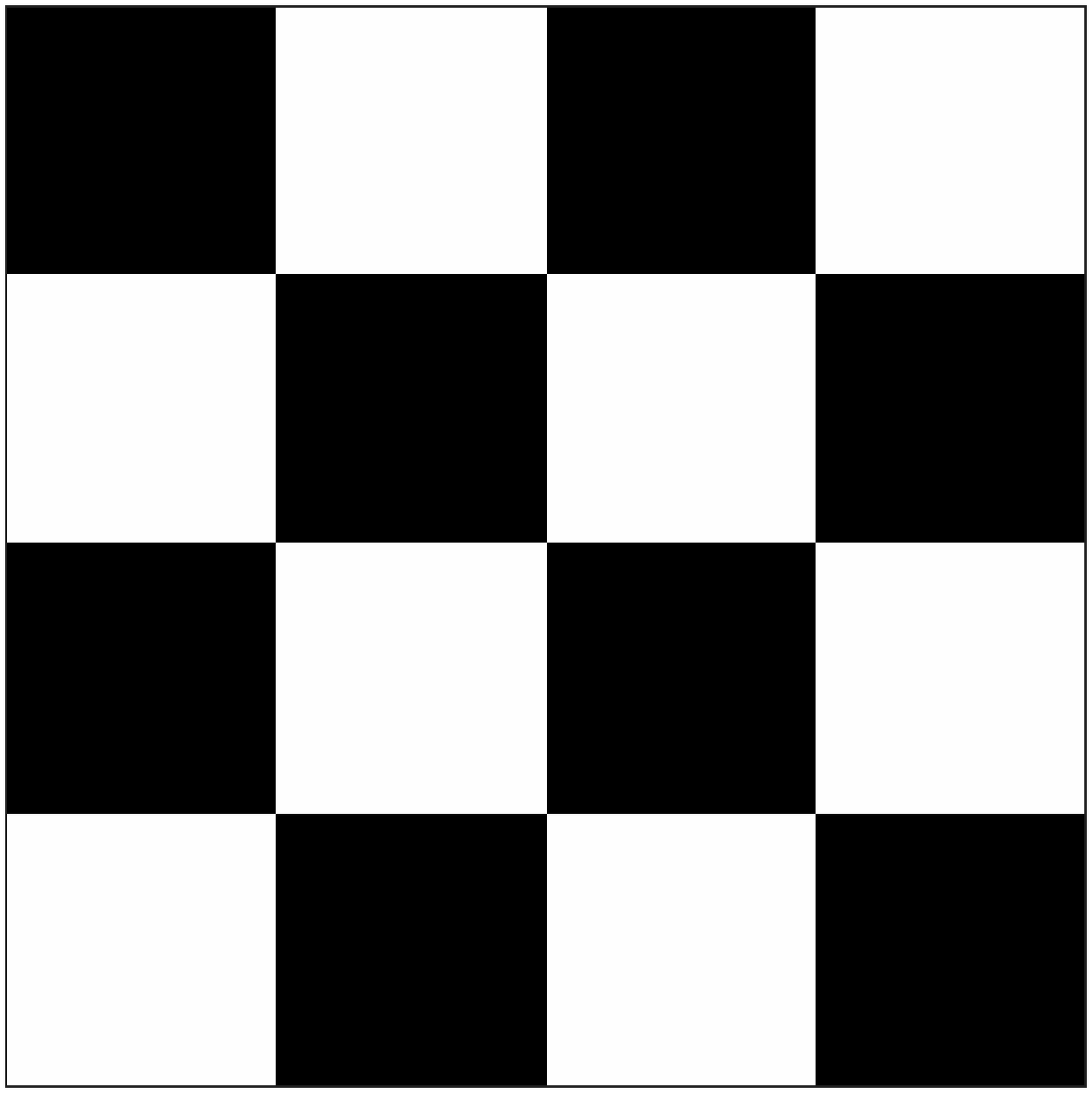} & 
\epsfxsize=4cm
\epsfbox{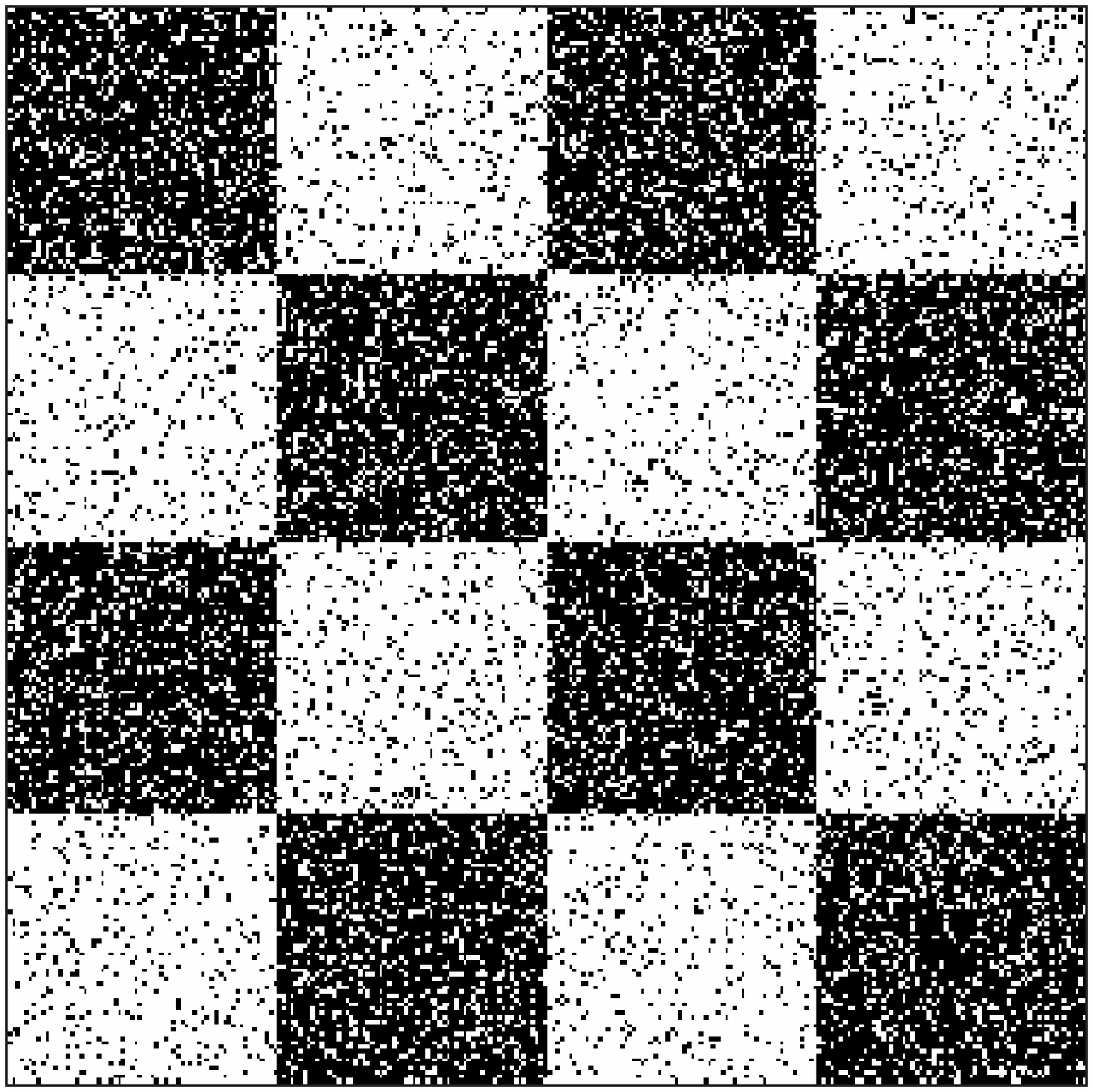}\\
\epsfxsize=4cm
\epsfbox{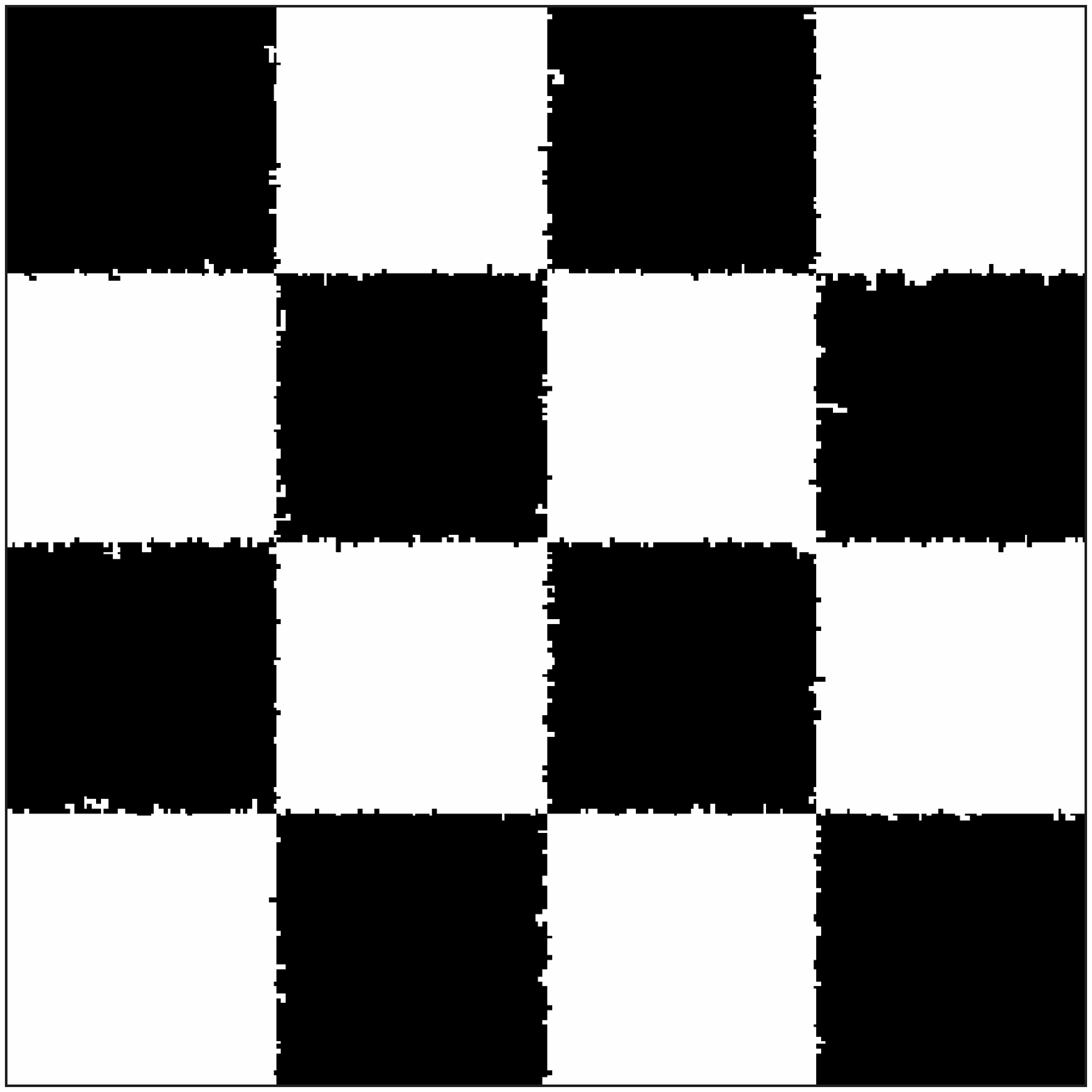} &
\epsfxsize=4cm
\epsfbox{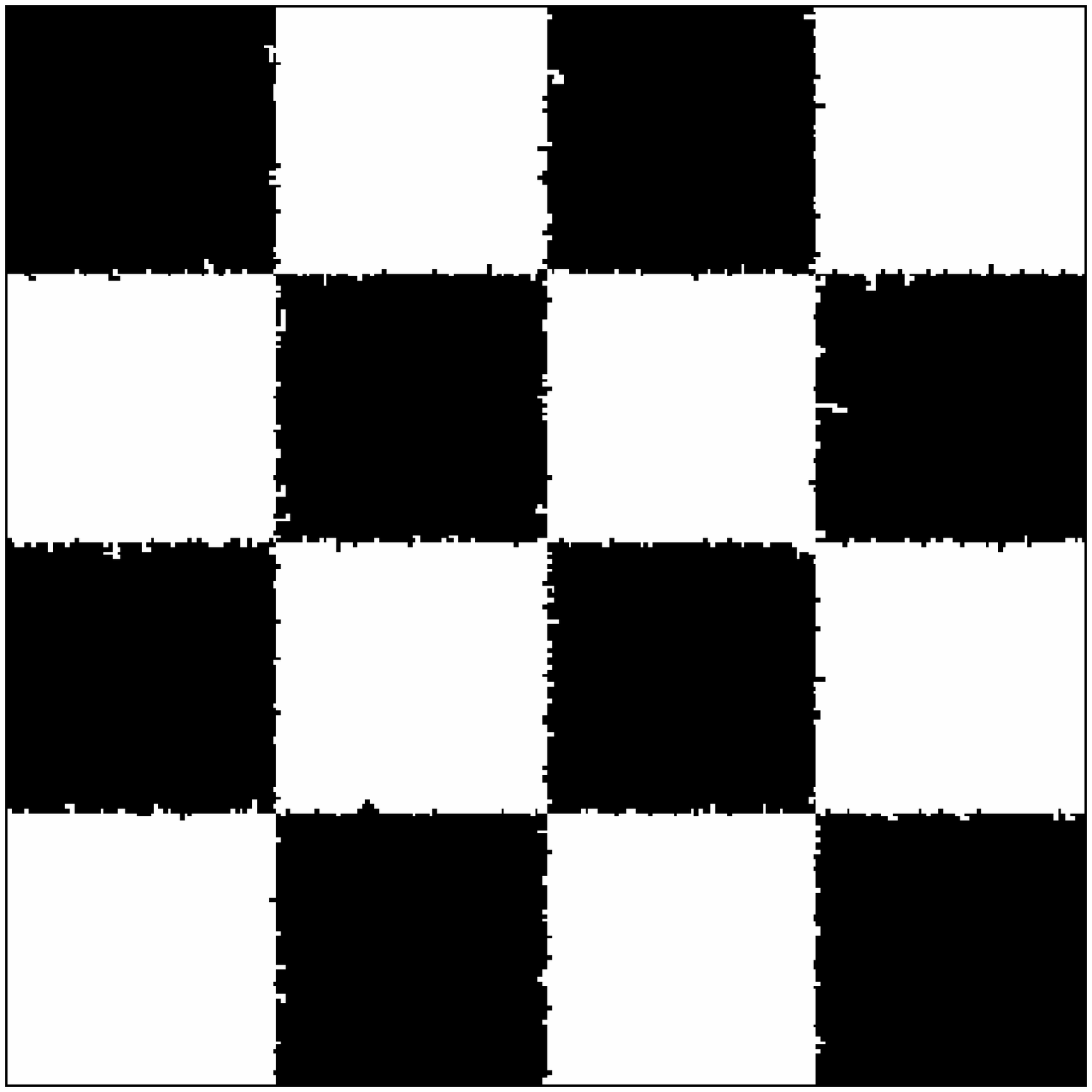} \\
\epsfxsize=5cm
\epsfbox{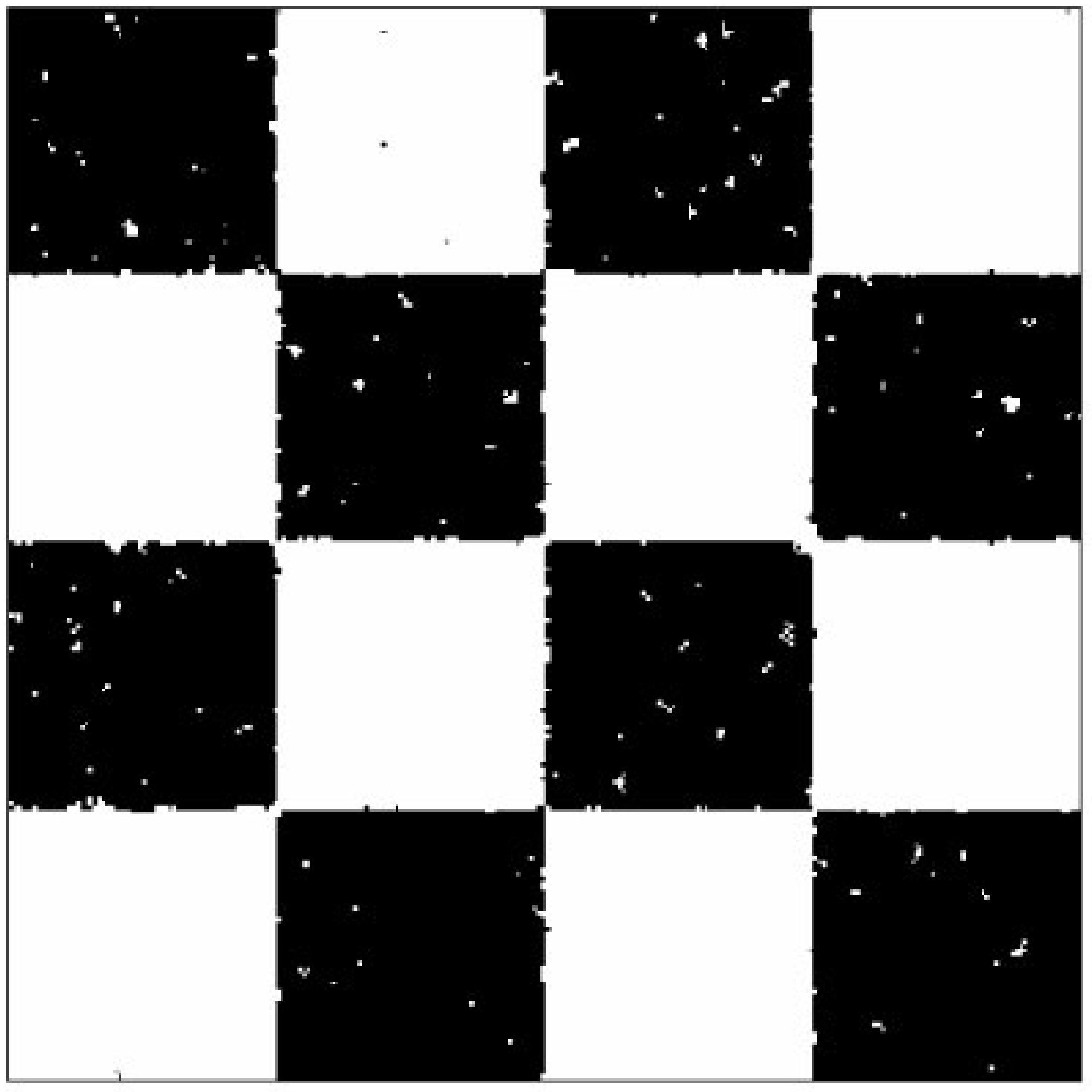} &
\epsfxsize=5cm
\epsfbox{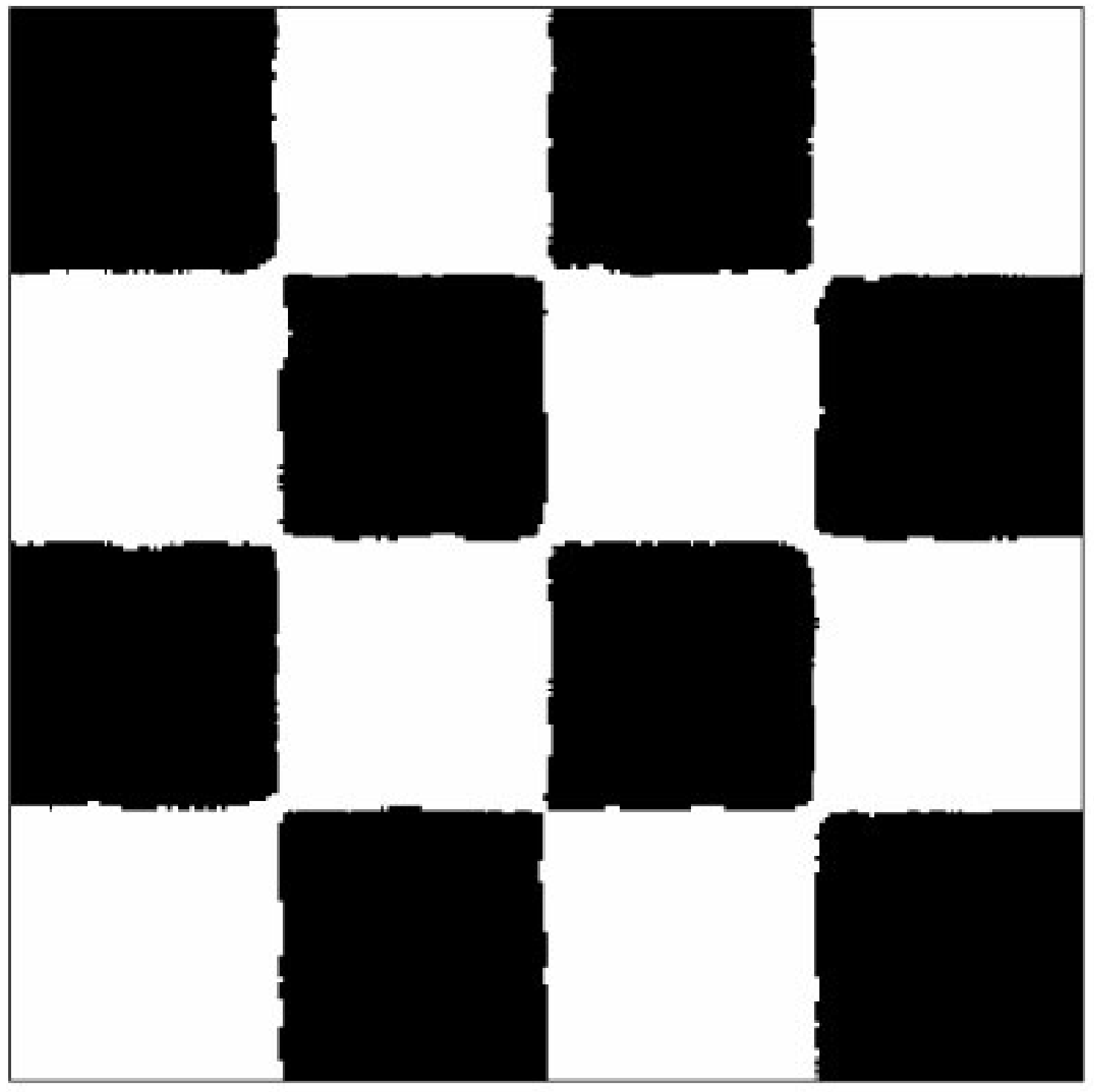}
\end{tabular}
\end{center}
\caption{{\it \small First row: on the left, the original binary image $I_0$ of size
  $256\times 256$; on the right, the corrupted image $I$. White pixels have been
  changed with probability $p=0.1$ and black ones with probability $q=0.2$.
  Second row: on the left, the result of the denoising algorithm with
  $\eps=10^{-2}$ when removing
  first black components and then white ones (i.e. applying $T^+_{s(n,q,\eps)}
  \circ T^-_{s(n,p,\eps)}$); on the right, denoising by first removing
  white components and then black ones (i.e. applying $T^-_{s(n,p,\eps)}
  \circ T^+_{s(n,q,\eps)}$). The two images are not the same, illustrating the
  fact that the two operators  $T^-_s$ and $T^+_s$ do not commute. However in
  both cases, small black and white components due to noise have been
  removed. Only the boundaries of the remaining ones are different.
  Third row: results obtained when applying median filtering with a disk of
  radius $r=2$ (on the left) and a disk of radius $r=5$ (on the right). The
  value of the parameter $r=2$ seems too small since some noise is still
  present in the black components. Now, for the value $r=5$, one can notice
  that the black corners have been eroded. The reason for this is that, in the
  noisy image, the probability parameters $p$ and $q$ (used to corrupt
  respectively the white and black pixels) were such that $q>p$.}} 
  \label{damier.fig}
\end{figure}

Let $I_0$ be the original (non degraded) binary image of size $n\times
n$. 
This original image $I_0$ is then corrupted by  impulse noise, which has a
probability parameter $p$ in the white components and a probability parameter
$q$ in the black 
ones (see Figure \ref{damier.fig} for an example). We shall see in the next
section why it is important to allow  
black and white pixels to be destroyed with a different probability. Thus the
noisy image $I$ is given by
\begin{equation}
\forall x , \hspace{0.2cm} I(x)=I_0(x)\cdot (1- \zeta_p(x))+ (1-I_0(x)) \cdot
\zeta_q(x) ,
\label{noise_binaire.eq}
\end{equation}
where the $\zeta_p(x)$'s (resp. $\zeta_q(x)$'s) are independent Bernoulli random variables
with parameter $p$ (resp. $q$). In other words, we have the following
conditional probabilities
$${\mathbb{P}}(I(x)=0 \, | \, I_0(x)=1)=p \hspace{0.4cm} \mathrm{and} \hspace{0.4cm}
{\mathbb{P}}(I(x)=1 \, | \, I_0(x)=0)=q .$$

\begin{figure}[!htbp]
\begin{center}
\epsfxsize=10cm
\epsfbox{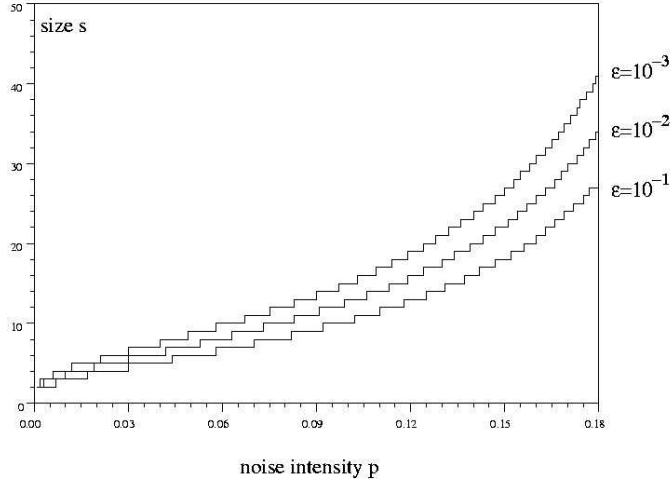}
\end{center}
\caption{{\it \small The size threshold $s(n,p,\eps)$ as a function of the 
noise probability parameter $p\in[0,0.18]$, for $n=256$ and $\eps=10^{-1}$, $10^{-2}$
and $10^{-3}$. }}
\label{seuils.fig}
\end{figure}

As can be seen in Figure \ref{damier.fig}, the impulse noise creates
small black and white connected components. These small components 
will be removed using a statistical decision based on their size 
(``size'',  in this paper, always means ``area'').
We are first interested in the black connected components (with respect to
4-connectivity). The results of the previous section give us 
the threshold function and also the probability of appearance of such
components. More precisely, the threshold function for a 
given (fixed) black component of
size $k$ is $\theta(n)=n^{-2/k}$ and  its asymptotic
appearance probability in a $n\times n$ image of noise with
probability parameter $p(n)=c\theta(n)$, as $n$ goes to infinity, is equal to 
$$1-e^{-c^k} . $$  
Now, if we are interested in the appearance of a component of size
$k$ (i.e. any of them, not only a given one), Theorem \ref{th:psi} claims that
the asymptotic (for large $n$) probability of appearance  is 
$$
1-e^{-a_k c^k} ~, 
$$
where $a_k$ is the number of 4-connected components one can make with exactly
$k$ pixels (up to translations). Writing $c^k=n^2 p^k$, 
we thus have an approximation for the probability
of appearance of a component of size $k$ in the $n\times n$ image,
with a proportion $p$ of black pixels. We denote by
$\mathrm{PA}(n,k,p)$ this approximation~:
$$
\mathrm{PA}(n,k,p) = 1-e^{-n^2 a_k p^k} ~.
$$ 
The $4$-connected components are known in the combinatorics
literature as
``square lattice animals'' or ``polyominoes''. 
Counting these objects is a difficult combinatorial problem
and there is no general expression for $a_{k}$. However, some 
asymptotic results are known: 
a concatenation argument \cite{klarner67} shows that there
exists a  constant $a$, called \textit{growth constant}, such that:
$$
\lim_{k \to \infty} (a_{k})^{\frac{1}{k}} = \sup_{k \geq 1}
(a_{k})^{\frac{1}{k}} = a ~.
$$ 
The exact value of $a$ is unknown. Numerical estimates give 
$a\simeq 4.06$ and the best published
rigorous bounds for it are $3.9 < a < 4.65$ 
(see \cite{conway95,jensen00,klarner73}). But thanks to
some numerical studies\footnote{for up-to-date information on the topic, see
  the web-site of the ``On-line Encyclopedia 
  of Integer Sequences'', {\tt http://www.research.att.com/$\sim$njas/sequences/} and
  references therein.}, 
the values of the sequence $(a_k)_{k\geq 1}$ are known up to
$k=47$, which will be enough in practice for denoising applications. 
The first terms are: $a_{1}=1$, $a_{2}=2$, $a_{3}=6$, $a_{4}=19$, etc. 
Furthermore, numerical computations show that for $p\leq p_{max}\simeq 0.2$,
one has $a_{k+1} p \leq  a_{k}$ for $k\in[1,47]$, 
which ensures that $\mathrm{PA}(n,k+1,p)\leq
\mathrm{PA}(n,k,p)$. 
This means that the probability of appearance of an animal is a
decreasing function of its size. This is rather reasonable: for fixed values
of $p$ and $n$, it would not make much sense to keep a  
connected component of size $k$ and to remove one of size $k'> k$.
\vskip 3mm
Let us fix a (small) positive real $\eps$ which will be our risk
probability, in the sense of statistical testing. If the size $k$ of a
connected component observed in a noisy image $I$
is such that $\mathrm{PA}(n,k,p)\leq \eps$, 
then we will consider that it comes from the original image $I_0$, 
and keep it. If 
$\mathrm{PA}(n,k,p)>\eps$, it will be regarded as noise and removed.  
Thus the size threshold for the components we keep is defined by:
\begin{equation}
s(n,p,\eps)=\mathrm{inf} \{ k \, ;\; 
\mathrm{PA}(n,k,p) = 1-e^{-n^2 a_k p^k}
\leq \eps \} .
\label{formule_seuil.eq}
\end{equation}
A component with size $k\geq s(n,p,\eps)$ has a very low probability
(less than $\eps$) of appearing in a pure noise image. 
In Figure \ref{seuils.fig}.b, we plot the size threshold
$s(n,p,\eps)$ as a function of the noise probability parameter $p$, for a fixed
value of $n=256$ and three different values of $\eps$: $10^{-1}$, $10^{-2}$
and $10^{-3}$.

\noindent 
The algorithm for the binary image denoising can be decomposed in four steps: 
\begin{enumerate}
\item  Extract all the 4-connected black components of the noisy image $I$.
\item  Remove the ones which have a size less than $s(n,p,\eps)$ (i.e. change
  their pixels into white). Obtain a new binary image.
\item Extract all the white 4-connected components of this new image.
\item  Remove the ones which have a size less than $s(n,q,\eps)$ (i.e. change
their pixels into black), to obtain the final denoised image denoted by
$\tilde{I}= T I$. 
\end{enumerate}

\noindent To summarize,  this denoising
filter $T$ can be written as: 
$$ T = T^+_{s(n,q,\eps)} \circ T^-_{s(n,p,\eps)} ~,$$
where  $T^+_s$ (resp. $T^-_s$) is the morphological area opening
(resp. closing) of size $s$ defined by L. Vincent in \cite{vincent}. 
See Figure \ref{damier.fig} for an example
of the obtained result and for a comparison with the results obtained with a
more standard binary filter (namely the median filter). 

Before explaining how
this method will be extended to grey level images, let us 
make a few general comments.

\bigskip

\begin{itemize}
\item The method is valid when $p$ is not too large, since we need $a_kp^k$ to
  be small. In practice, we are limited to $p\leq p_{max}\simeq 0.2$.
\item The dependence on $\eps$ is low since it is in fact a
  $\log(\eps)$-dependence.  Indeed, $1-e^{-n^2 a_k p^k}$ is approximately
  equal to $n^2 a_k
p^k$ when the value of this expression is small. If we replace $a_k$
  by $a^k$,  the threshold for
  the minimal size of the components 
  we keep is approximately given by
$$ 
s(n,p,\eps)\simeq \frac{\log\eps -2\log n}{\log a +\log p} ~.
$$
\item The boundaries of the remaining components are not smoothed. This 
comes from the fact that when some noise is at the boundary of a component, it
becomes part of it.  
In order to remove it, one would need other {\it a priori} knowledge of the
original image (such as smooth or  
straight boundaries as in the case of Figure \ref{damier.fig}).  
It is actually a general problem of image denoising: one has to define some
{\it a priori} model for the image.  
Here the underlying model is that the original binary image is made of
``large'' (as compared to the noise) black and white 
connected components.  
\item The two filters $T^-_s$ and $T^+_s$ do not commute (see Figure \ref{damier.fig}). This was already
  noted by  Vincent in \cite{vincent}.  
One solution he proposed is to use them in alternating sequential filters \cite{Serra2, vincent2} 
with increasing sizes of area. This may not be a real issue, 
since he also noticed 
that $T^-_s \circ T^+_s$ and $T^+_s \circ T^-_s$ are visually 
extremely close (this will be even more true for grey level images). 
Another solution, proposed by Masnou and Morel in \cite{masnou98}, and
  then formalized by Monasse \cite{monassePhD}, is to process simultaneously
  upper and lower level sets. This grain filter denoted $G_t$ (where $t$ is
  the area threshold) is done by a pruning of the tree of all level
  sets, built thanks to the inclusion principle (this algorithm which is very
  fast is called the Fast Level Set Transform \cite{monasse}). \\
 The fact that the
foreground and the background of a binary image are treated in a
complementary way is a general problem in Mathematical Morphology. Many
  operators are not self-dual, and they often occur pairwise: like dilation/erosion
  and opening/closing
for example. It is worth mentioning that in \cite{Heijmans2}, H. Heijmans describes a
general method to construct morphological operators which
are self-dual.
\item It is generally considered that, for consistency reasons, using
the 4-connectivity
on the black (or white) pixels should be
followed by using the 8-connectivity for the
complementary set. From a theoretical point of view, the method we proposed
can be extended to 8-connectivity in a straightforward way. 
In order to apply the method, one would have to count the number of
8-connected components of size $k$, which are not available in the literature,
whereas the $(a_{k})$'s are known up to $k=47$.
Consequently, for our application to image denoising, we decided to treat the
foreground and 
the background in the same way, with 4-connectivity.
\item If the original image $I_0$ is all white ($I_0\equiv 0$), and if it is
  corrupted by some noise with probability parameter $p$  
as described by equation (\ref{noise_binaire.eq}), 
we obtain an image $I$ which is ``pure noise''. The probability  
that it contains a connected component with size larger than $s(n,p,\eps)$ is
(by definition of $s(n,p,\eps)$ and thanks  
to theorem \ref{th:psi})  less than $\eps$. Thus, 
$$
{\mathbb{P}}(T^-_{s(n,p,\eps)} I = I_0) \geq 1-\eps ~,
$$
which means that, with probability larger than $1-\eps$, pure noise is
completely removed.  \\
Thus $\eps$ represents the ``significance level'' of our statistical method
for denoising: the probability of not removing a component coming from the
noise is less than $\eps$. In practice, we generally take $\eps=10^{-3}$ or
$10^{-2}$ (the results are visually the same). This 
parameter $\eps$ is completely independent of the image (which is not the case
of the size $n$ of the image, or the probability parameter $p$ of the impulse
noise): $\eps$ has to be fixed by the user in the same way as the risk level in
statistical hypothesis testing.
\item For reasons of simplicity, the denoising filter $T$ was described in the
  framework of an image of
  size $n\times n$. The method extends straightforwardly to an image of size
  $m\times n$, where $m\neq n$, by simply changing the size threshold
  $s(n,p,\eps)$ into $\tilde{s}(m,n,p,\eps)$ defined by
$$\tilde{s}(m,n,p,\eps)=\mathrm{inf} \{ k \, ;\; 
\mathrm{PA}(\sqrt{nm},k,p) = 1-e^{-m n a_k p^k}
\leq \eps \} .$$

\end{itemize}

\subsection{Grey level images}
\label{part_grey.subsec}

Let $u$ be a grey level image, of  size $n\times n$ and grey level values in
the range $[0,255]$. Assume that this image is corrupted by impulse noise
with probability parameter $p$. This means that the observed noisy image $v$ may be
written in the form:
\begin{equation}
\forall x , \hspace{0.2cm} v(x)=(1-\zeta_p(x)) \cdot u(x) + \zeta_p(x) \cdot \nu(x) ,
\label{noise_grey.eq}
\end{equation}
where the $\zeta_p(x)$'s are independent Bernoulli random variables with 
parameter $p$ and the $\nu(x)$'s are i.i.d.r.v.'s, uniformly distributed 
on $[0,255]$.

\begin{figure}[h]
\begin{center}
\begin{tabular}{ccc}
\epsfxsize=5.6cm
\negthickspace \negthickspace \negthickspace \negthickspace
\negthickspace \negthickspace \negthickspace \negthickspace
\epsfbox{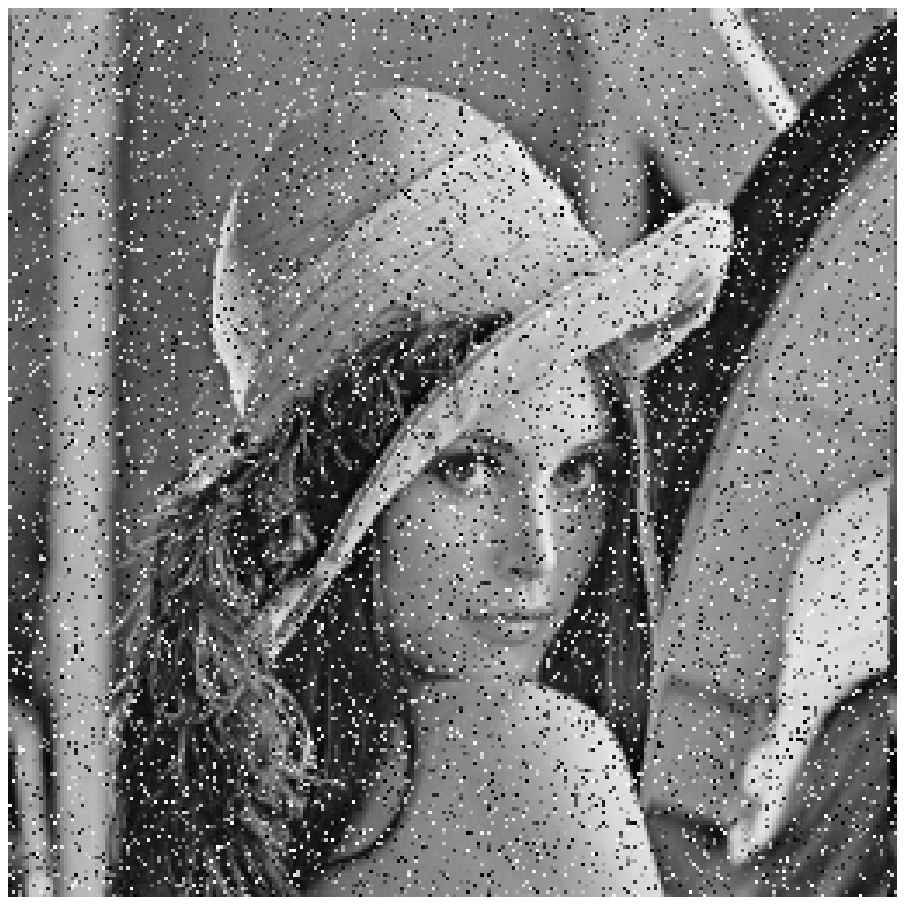} & 
\epsfxsize=5.6cm
\negthickspace
\epsfbox{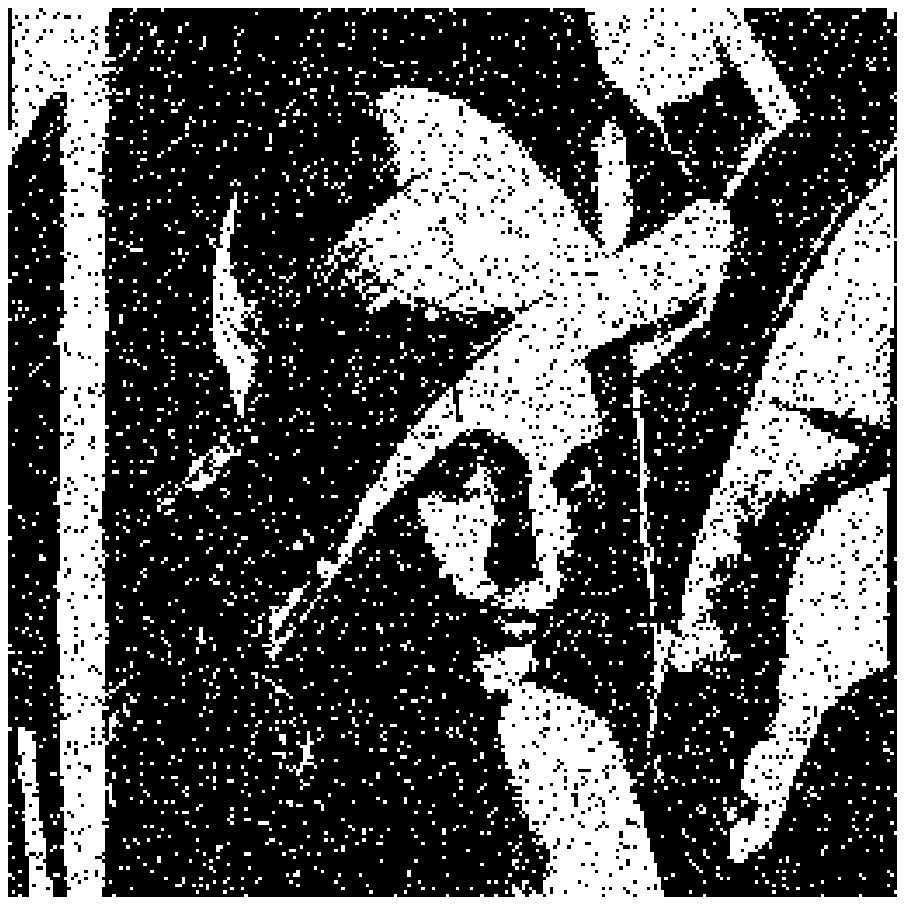} &
\epsfxsize=5.6cm
\negthickspace
\epsfbox{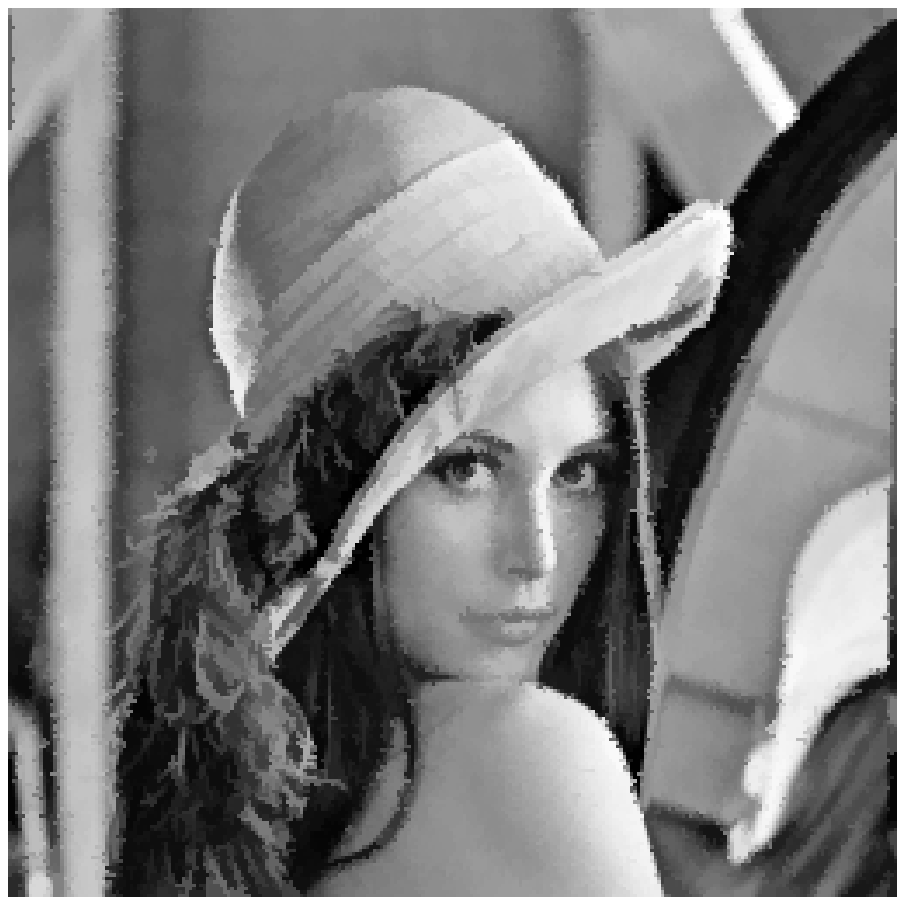}
\end{tabular}
\end{center}
\caption{{\it\small Left: image $v$ obtained with impulse noise with
  probability parameter $p=0.15$ on the Lena 
  image. Middle: thresholded image $v_\lambda$ for the grey level
  $\lambda=150$. Right: denoised image $\tilde{v}$ obtained by the
  noise-adapted grain filter $T$ with $\eps=10^{-3}$.}} 
\label{lena.fig}
\end{figure}

\begin{figure}[!htbp]
\begin{center}
\begin{tabular}{lr}
\epsfxsize=6cm
\epsfbox{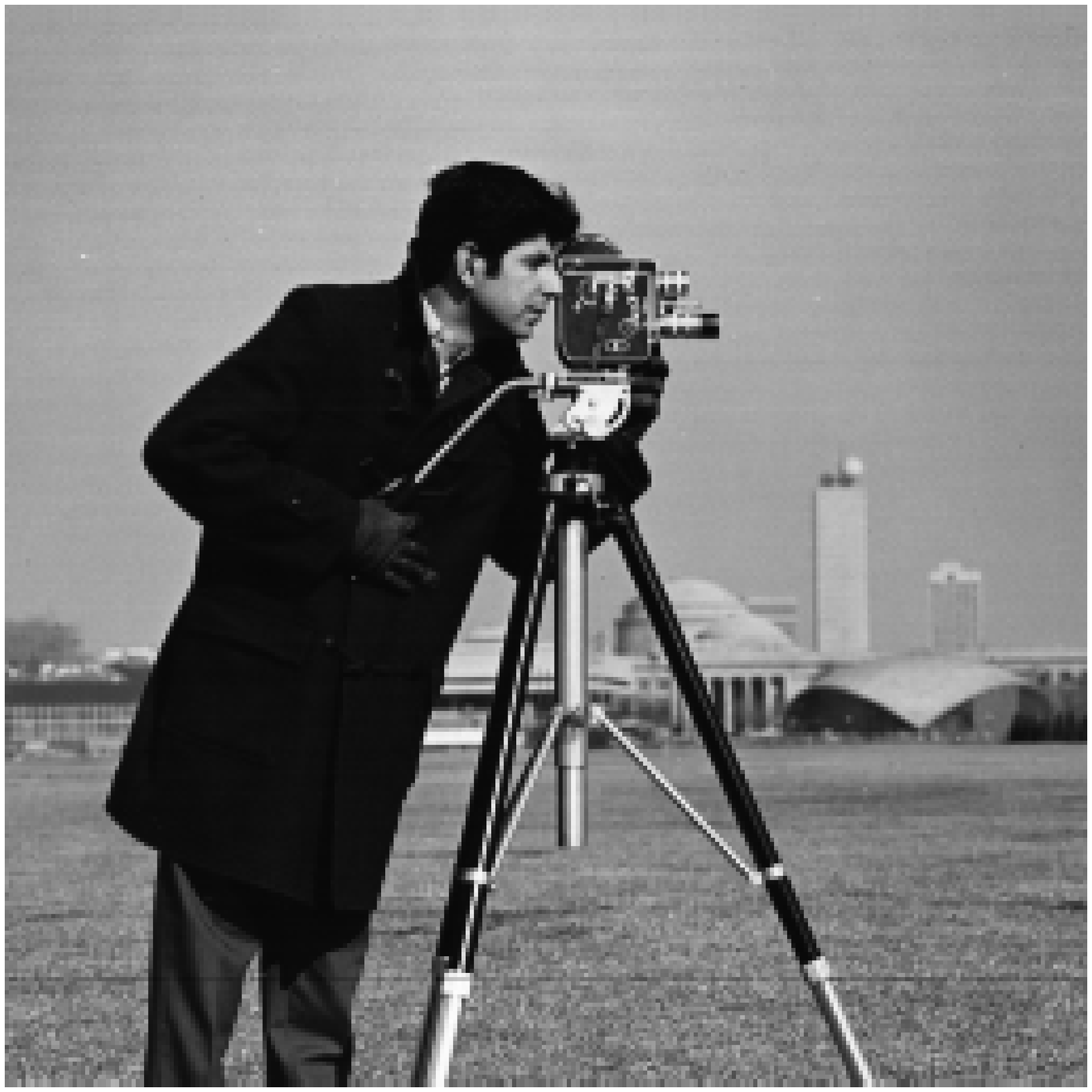} &
\epsfxsize=6cm
\epsfbox{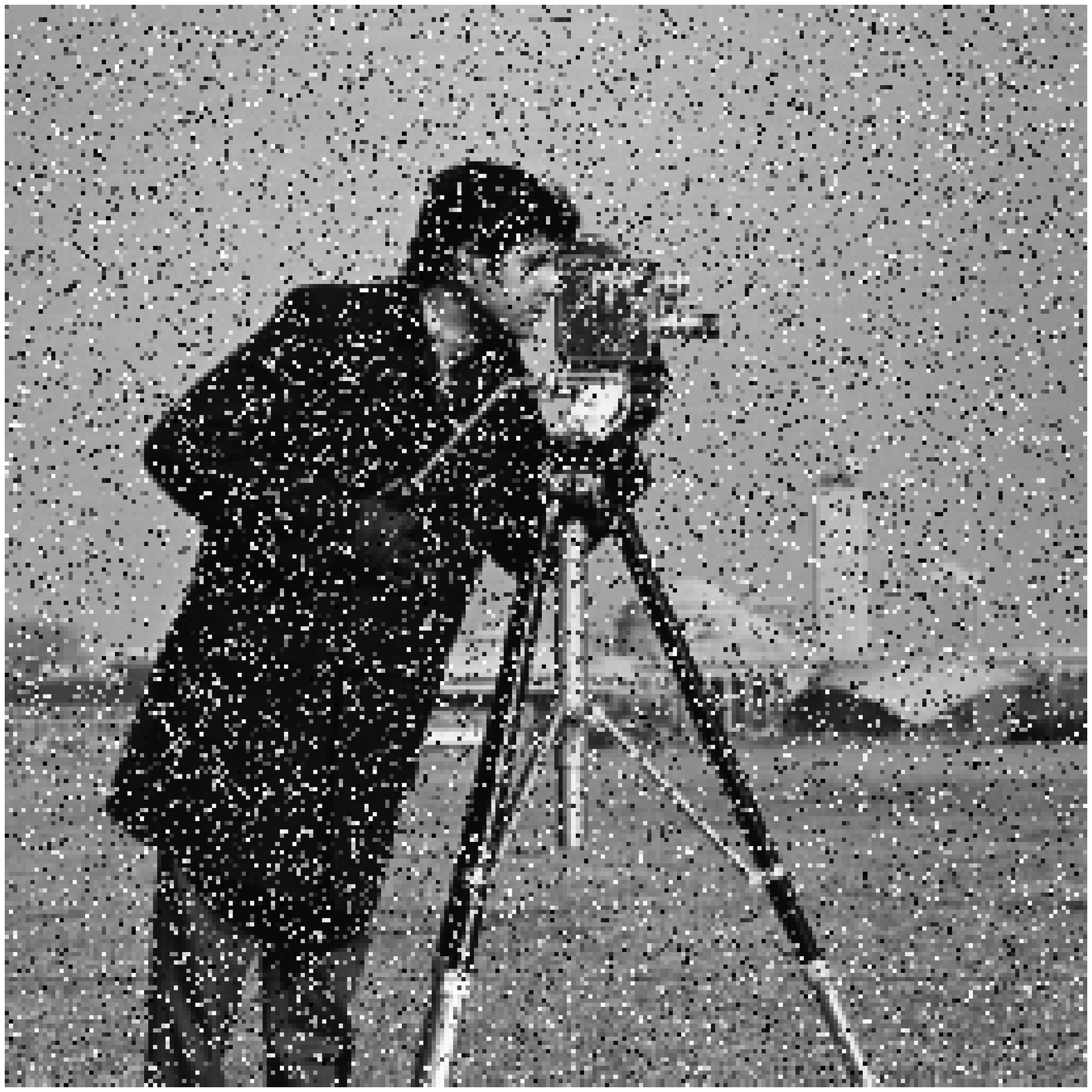}
\end{tabular}
\vskip 0.3cm
\begin{tabular}{lr}
\epsfxsize=6cm
\epsfbox{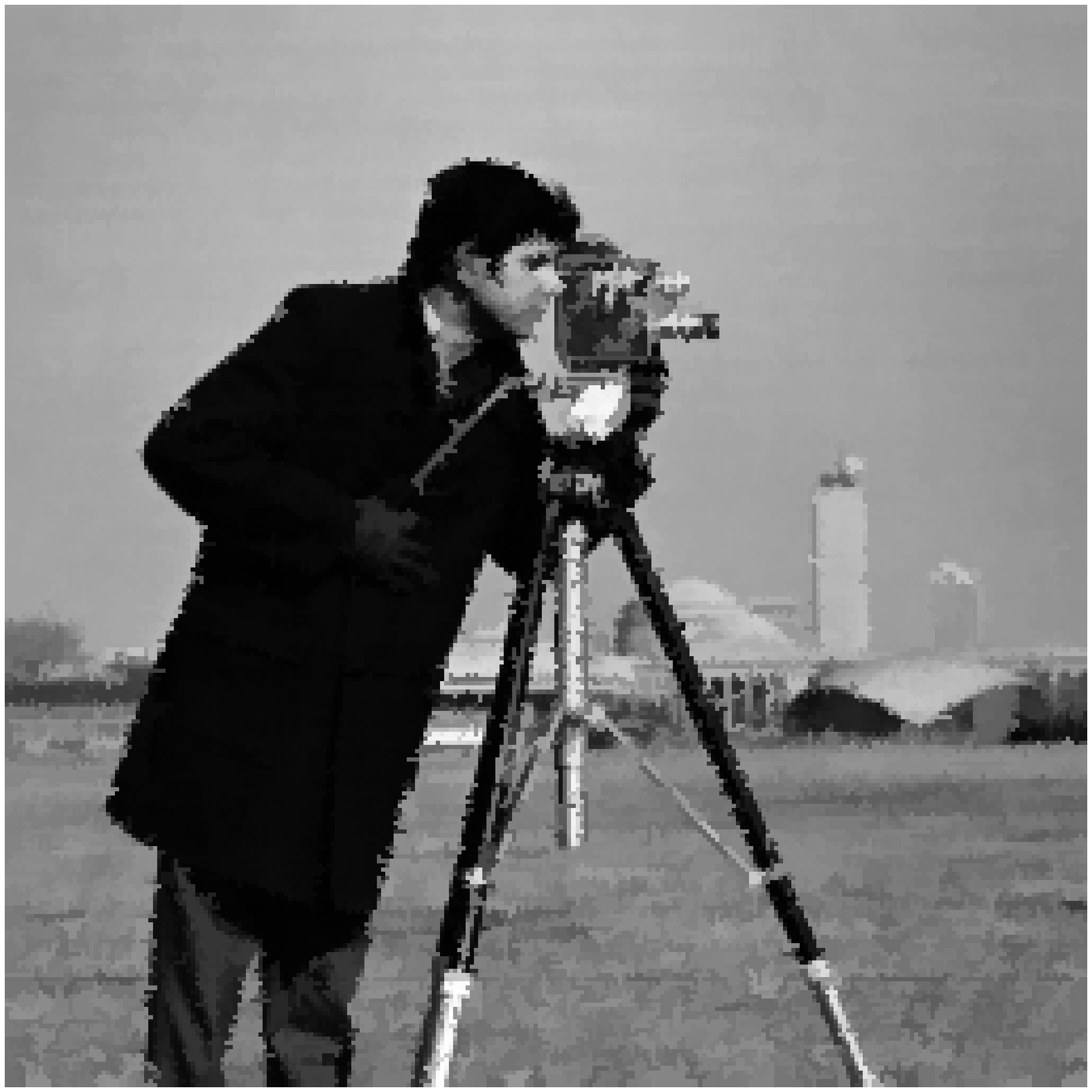} &
\epsfxsize=6cm
\epsfbox{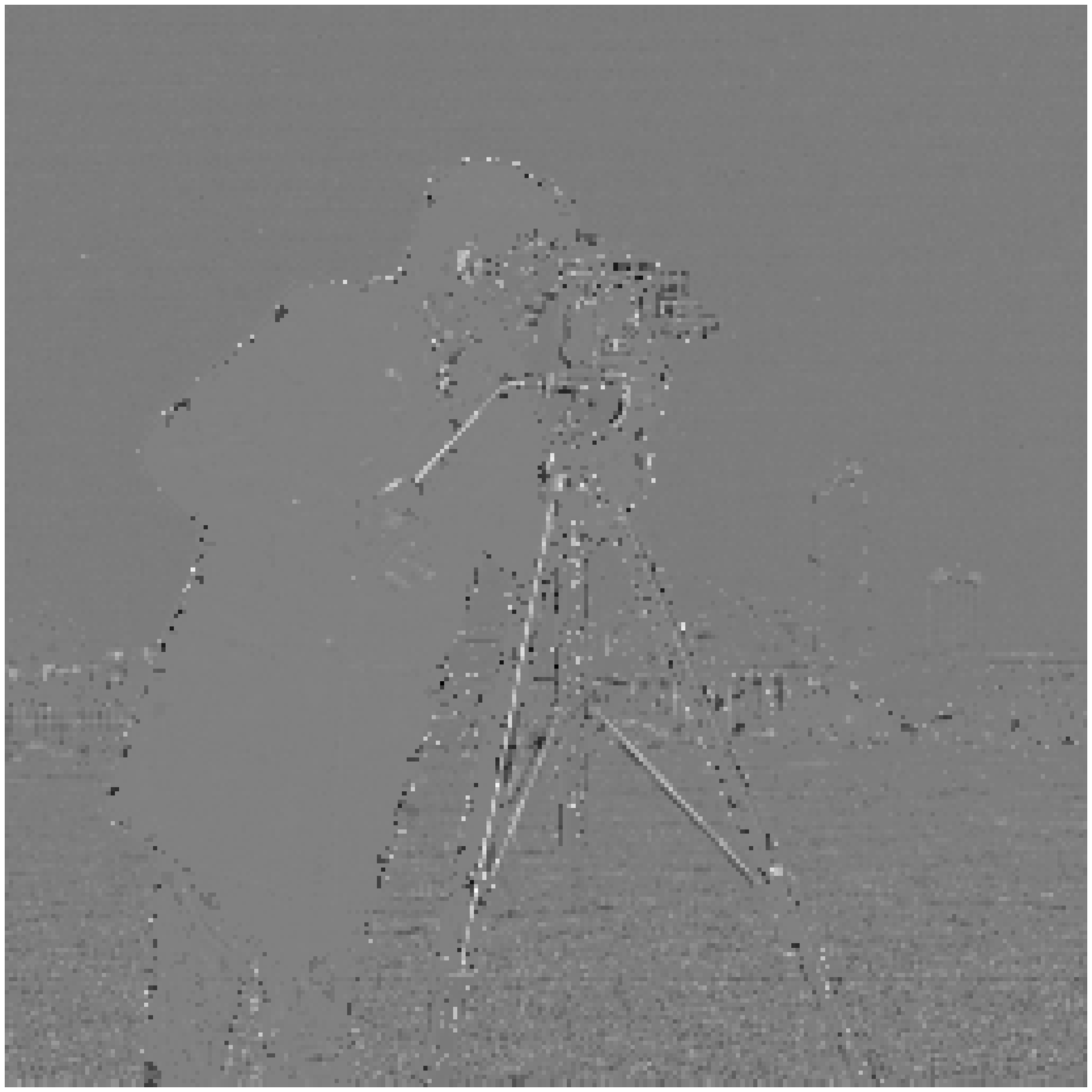}
\end{tabular}
\end{center}
\caption{{ \it \small From left to right, top to bottom: (a) the original cameraman image $u$
  (size $256\times 256$); (b) degraded image $v$, with impulse noise probability parameter
  $p=0.2$; (c) filtered image $\tilde{v}$, obtained with $\eps=10^{-3}$; (d)
  image of the difference $u-\tilde{v}$. It shows that most of the noise has
  been removed, except at the boundaries of the objects and also in the grass
  texture.}}
\label{cameraman.fig}
\end{figure}

For each level $\lambda\in[0,255]$, we can consider the thresholded images
$u_\lambda=\mathbb{I}_{ u\geq \lambda}$ and $v_\lambda=\mathbb{I}_{ v \geq
\lambda}$. The grey level images may then simply be recovered by
$u=\sum_\lambda u_\lambda$ and $v=\sum_\lambda v_\lambda$. 
The binary noisy image $v_\lambda$ is a corrupted
version of the binary image $u_\lambda$; they are related by  
$$
{\mathbb{P}}(v_\lambda (x)=0 \, | \,
u_\lambda(x)=1)=p\times\frac{\lambda}{256} 
\hspace{0.3cm} \mathrm{and}\hspace{0.3cm} 
{\mathbb{P}}(v_\lambda (x)=1 \, | \, u_\lambda
(x)=0)=p\times(1-\frac{\lambda}{256}) .$$
We are thus back in the framework described for binary images with parameters
$p_\lambda= p\lambda/256$ and $q_\lambda= p(1-\lambda/256)$. The image
$v_\lambda$ can be denoised following the method described in the previous subsection.
Finally, we reconstruct a grey level image by simply adding the binary
ones: $\tilde{v}=\sum_{\lambda} \tilde{v}_\lambda$.
This can be summarized by the formula
\begin{equation}
\tilde{v}= Tv = \sum_{\lambda=0}^{255} T^+_{s(n,q_{\lambda},\eps)} \circ \, 
T^-_{s(n,p_{\lambda},\eps)} (v_{\lambda} ) , \hspace{0.3cm} \mathrm{where} \hspace{0.2cm}
p_\lambda= p\frac{\lambda}{256} \hspace{0.2cm} \mathrm{and} \hspace{0.2cm} q_\lambda=
p\left( 1- \frac{\lambda}{256} \right) .
\label{eqfilter.eq}
\end{equation}
Figures \ref{lena.fig} and \ref{cameraman.fig} give two examples of
results obtained by this filtering. 

\medskip

One natural question that can be asked is whether the filter $T$ defined by 
formula (\ref{eqfilter.eq}) is a morphological filter (see \cite{Serra} and
\cite{vincent2} and references therein for the definition and properties of
morphological filters).
Unfortunately, the answer is negative. For two grey levels $\lambda\geq
\lambda'$, one has ${v}_\lambda \leq {v}_{\lambda'}$, and for a fixed area
threshold $t$ one would have  $T^-_t({v}_\lambda) \leq T^-_t({v}_{\lambda'})$
(because area openings and closings are morphological operators). Now, the two
thresholds $s(n,p_{\lambda},\eps)$ and $s(n,p_{\lambda'},\eps)$ can be
different, i.e.   $s(n,p_{\lambda},\eps) > s(n,p_{\lambda'},\eps)$ and thus
it is not necessarily true that
$T^-_{s(n,p_{\lambda},\eps)} ({v}_\lambda) \leq
T^-_{s(n,p_{\lambda'},\eps)}({v}_{\lambda'})$. This happens when
${v}_{\lambda}$ and ${v}_{\lambda'}$ both contain the same small black connected
component of size $k$ such that $s(n,p_{\lambda},\eps) > k >
s(n,p_{\lambda'},\eps)$. However, in the experimental results, we noticed that
this rarely happens: for most values of $\lambda$, one has
$\tilde{v}_{\lambda} \leq \tilde{v}_{\lambda-1}$. 

\begin{figure}[!htbp]
\begin{center}
\begin{tabular}{lr}
\epsfxsize=6cm
\epsfbox{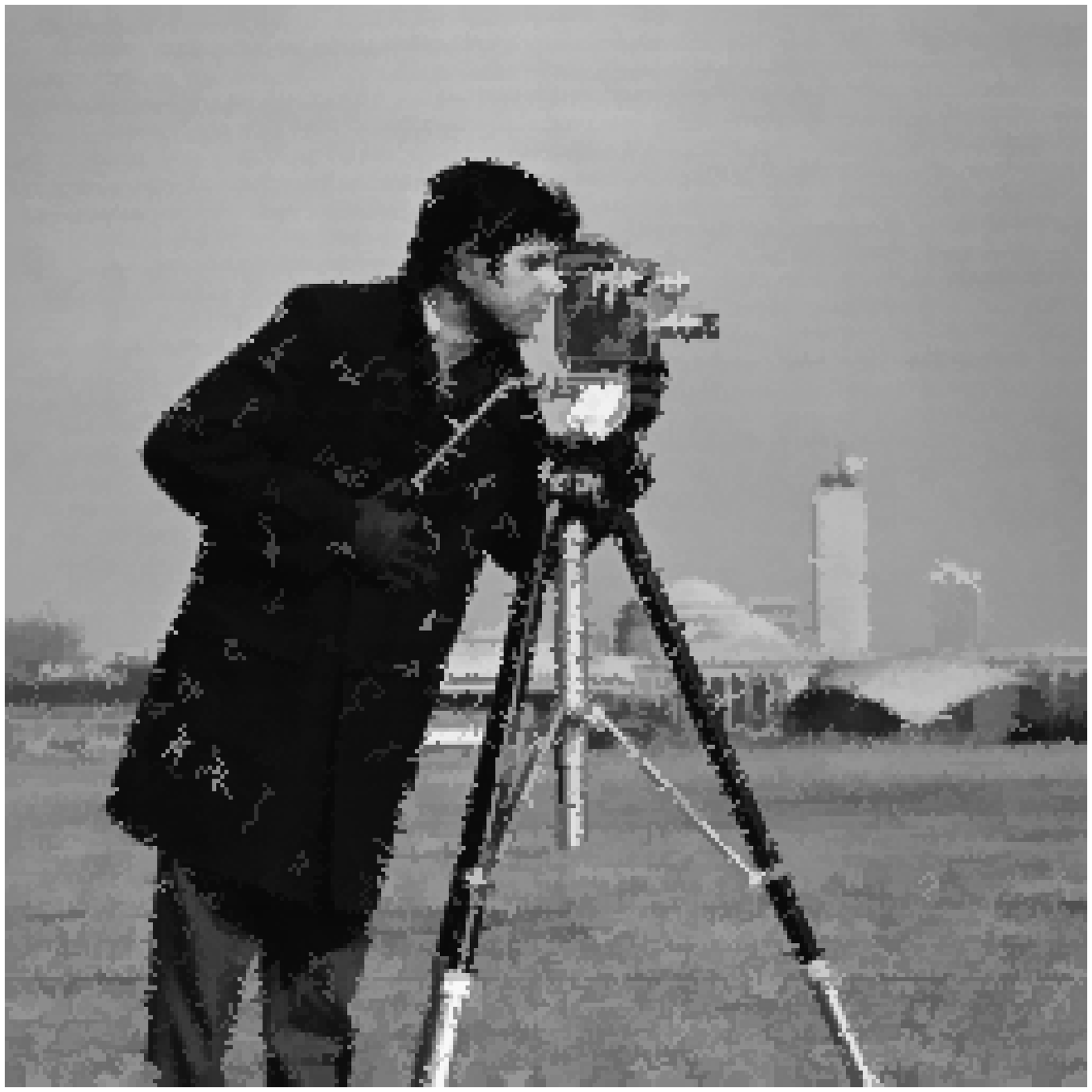} &
\epsfxsize=6cm
\epsfbox{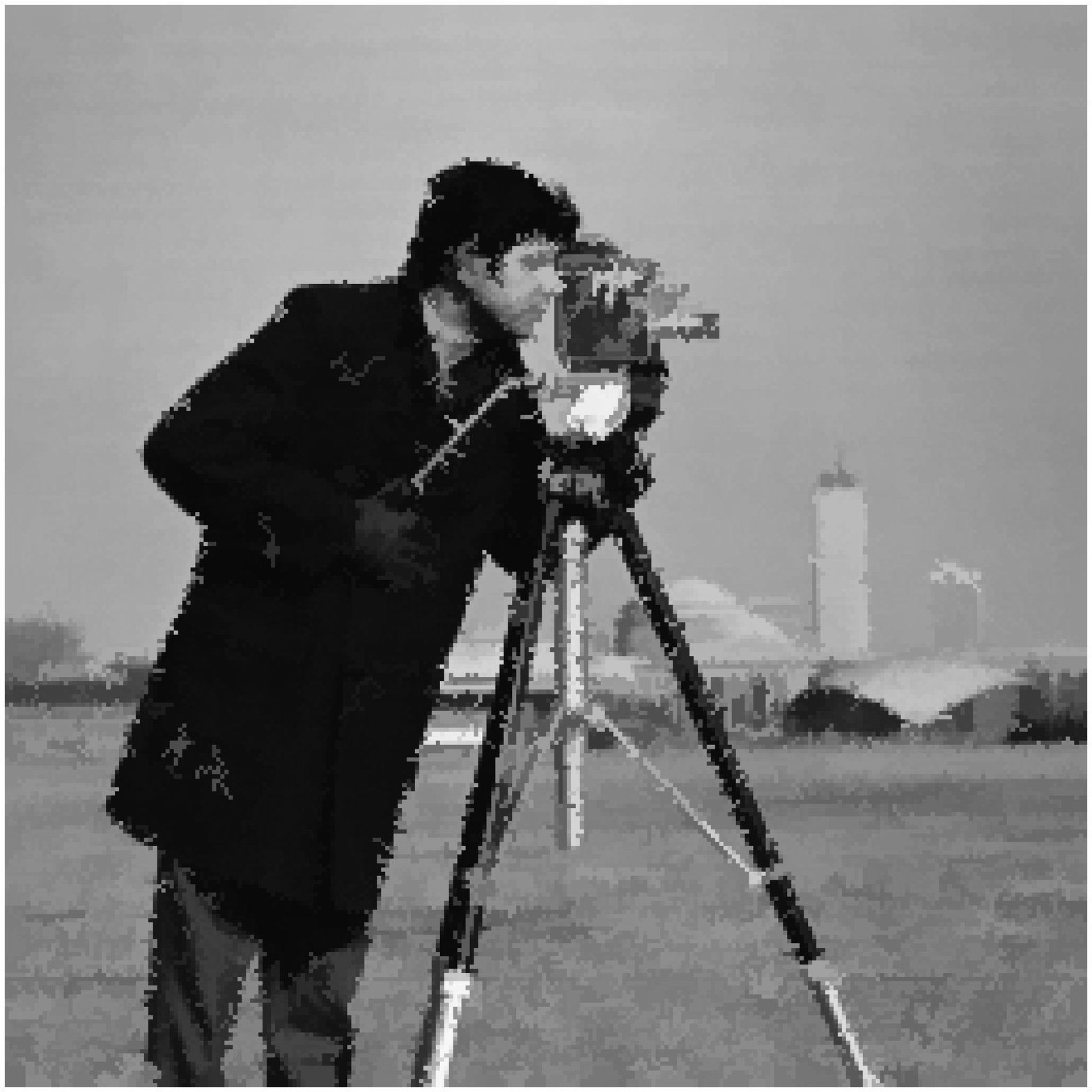}
\end{tabular}
\end{center}
\caption{{\it \small Result of the filtering of the noisy image $v$ with the
  usual grain filter $G_t$
  with area threshold $t=10$ on the left and $t=20$ on the right. }}
\label{compar_grain.fig}
\end{figure}

In order to illustrate the interest of adapting the area threshold to each grey level,
we treated 
the same image using our method, then using a fixed area threshold (for this
we used the algorithm developed by Monasse in \cite{monassePhD}). The results
are those of Figures \ref{cameraman.fig} and \ref{compar_grain.fig}.
Figure \ref{compar_grain.fig} shows the result of the usual grain filter,
denoted by $G_t$, for two different values of the area threshold: 
$t=10$ and $t=20$. 
One can notice that the parameter value  $t=10$ seems too
low since there is still some remaining noise (for example on the coat of the
cameraman). On the other hand the value $t=20$ seems too large, since some
of the original structures have disappeared (for instance 
the white parabola at the top of the building) and still too low (there is some
remaining noise on the coat).
These results have to be compared with the one of Figure
\ref{cameraman.fig}-c. This last figure shows that thanks to the adapted
area threshold $s(n,p_\lambda,\eps)$ a small white component can be kept and
at the same time, a larger grey component removed.  
These results also illustrate what we have proposed in this paper, namely an
adapted and automatic way to choose the 
right parameter for the area openings and closings. 

\subsection{Extension to other noise models}

\begin{figure}[!htbp]
\begin{center}
\begin{tabular}{ccc}
\hspace{-3cm}
\epsfxsize=7.5cm
\epsfbox{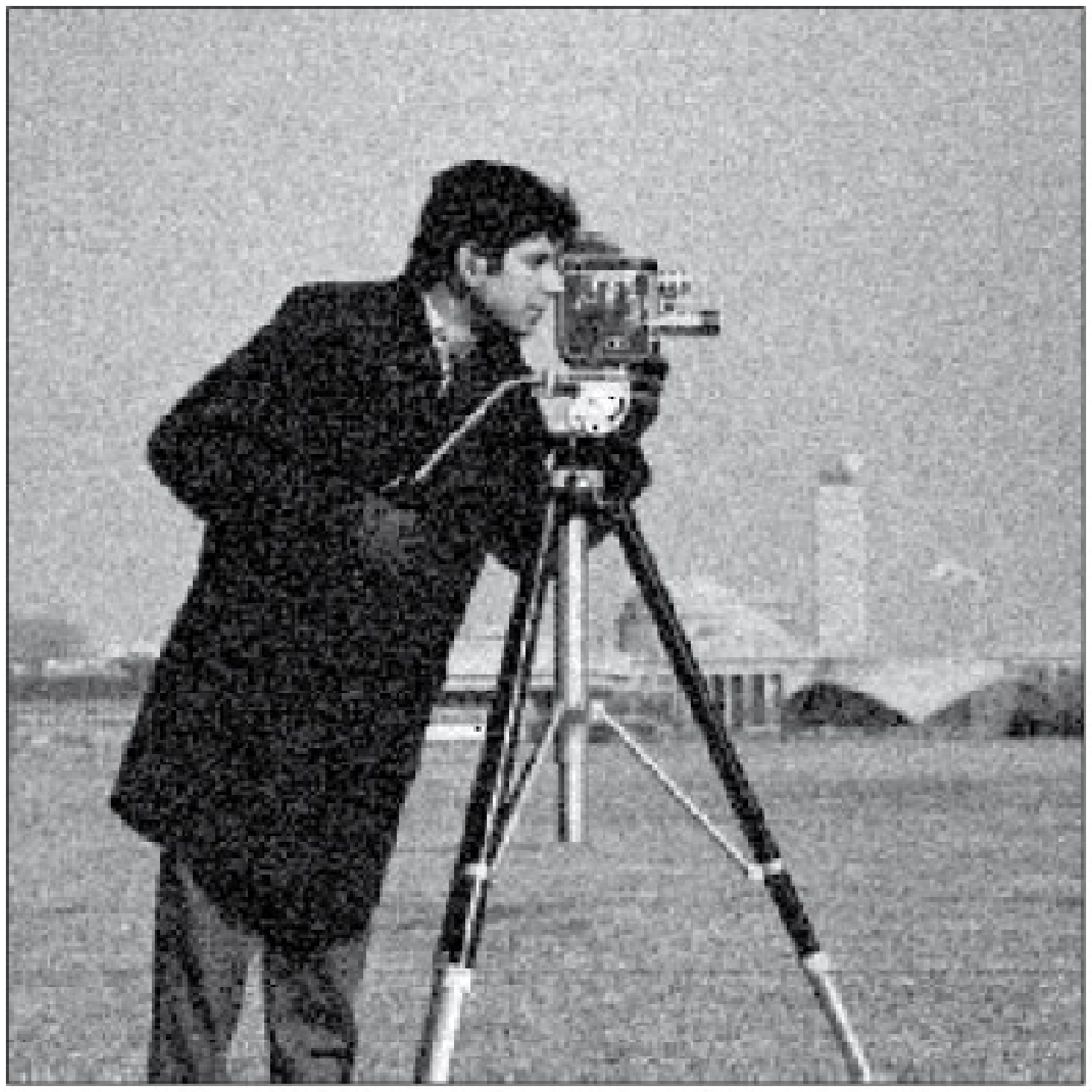} &
\hspace{-1cm}
\epsfxsize=7.5cm
\epsfbox{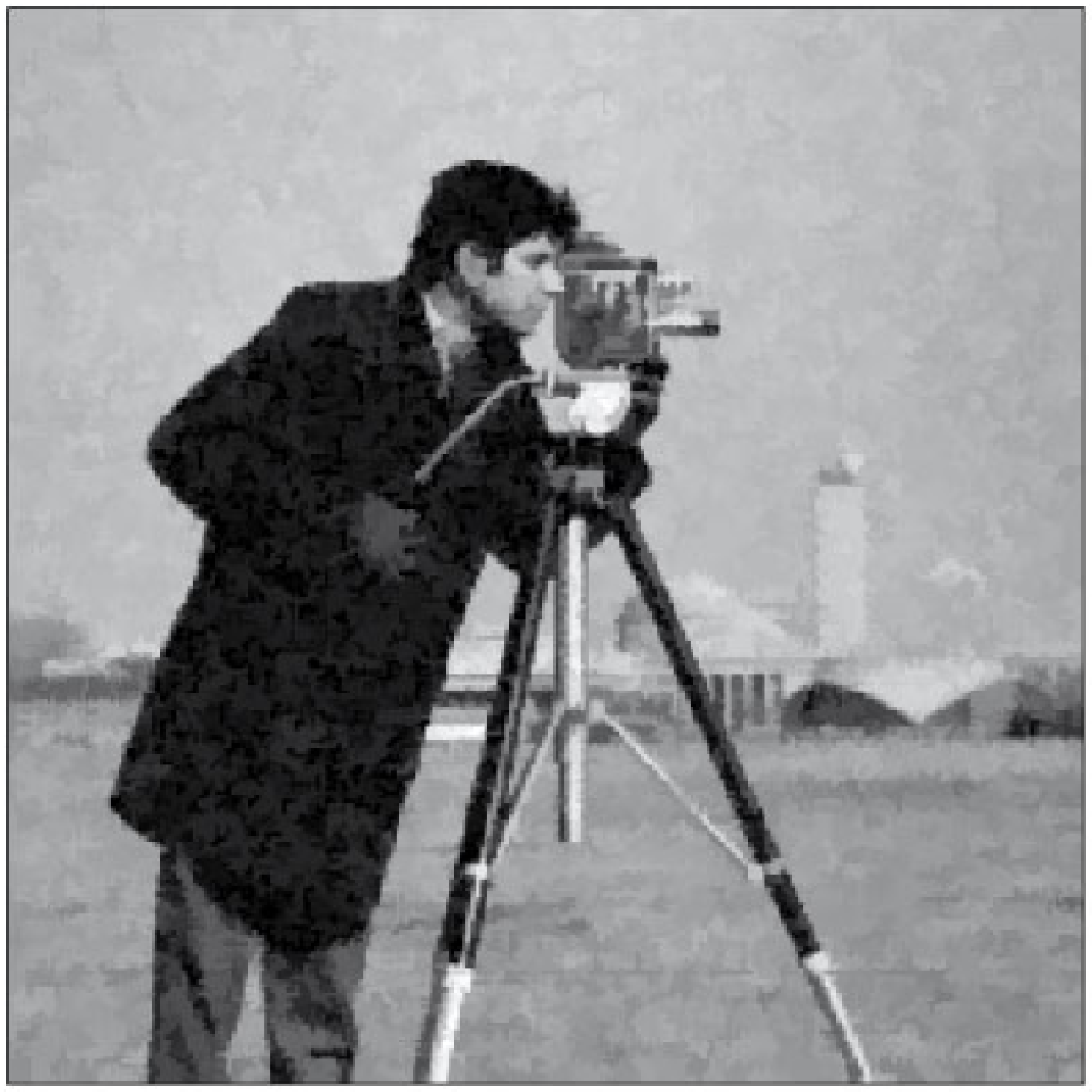} &
\hspace{-1cm}
\epsfxsize=7.5cm
\epsfbox{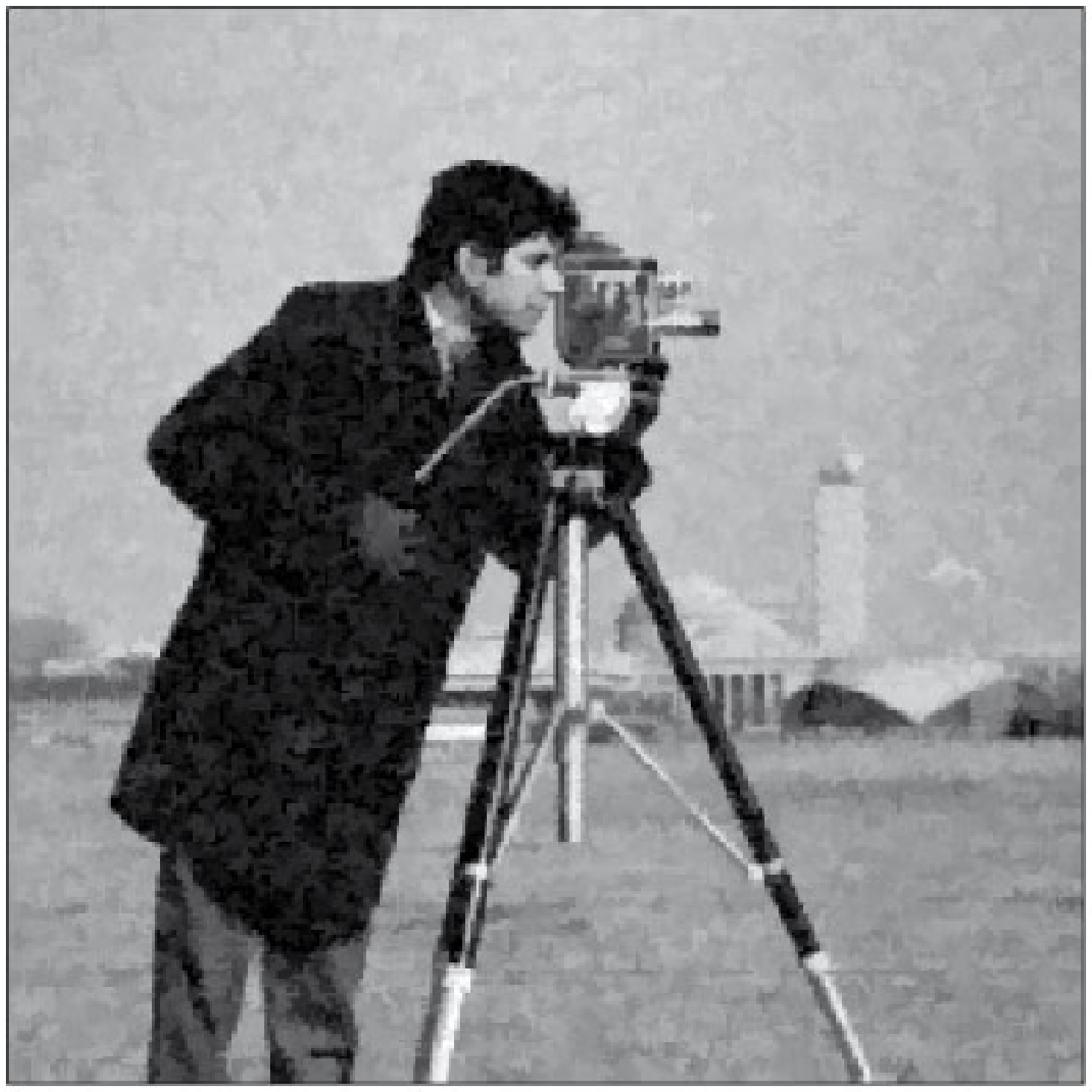}
\end{tabular}
\end{center}
\caption{ \it \small From left to right: (a) the cameraman image degraded by
     white noise with standard deviation $\sigma=15$; (b) denoising of the previous
    image by the filter defined by Equation (\ref{eqfilter.eq})
    with parameter values $p=0.2$ and $\eps=10^{-2}$; (c) denoising by the
    same filter with parameter values $p=0.1$ and $\eps=10^{-2}$.}
\label{camerbblanc.fig}
\end{figure}

In the previous subsection, we have explained how the theoretical results of
Section \ref{part_pa.sec} can be used to denoise an image degraded by impulse
noise. 
Now, even if the proposed denoising procedure corresponds to an impulse noise
model, it is interesting to see how it works in the presence of white noise.
An example of the obtained results is shown on Figure \ref{camerbblanc.fig}:
we used again the cameraman image, which is here degraded by
     white noise with standard deviation $\sigma=15$. It is then denoised
     using the filter defined by Equation (\ref{eqfilter.eq}) with parameter
     value $p=0.2$ and $\eps=10^{-2}$ (on the middle image) and $p=0.1$ and
     $\eps=10^{-2}$ (on the right image). The main question is here: how to
     choose the value of the parameter $p$ used in the filter, in relation to
     the standard deviation $\sigma$ of the white noise~?

In the case of impulse noise, we were able to relate the size threshold $s$ of
the area openings and closings to the impulse noise probability parameter $p$ and to
the grey level $\lambda$. The main result was then: if we take $u=0$ in
Equation (\ref{noise_grey.eq}), then the degraded image $v$ is pure impulse
noise, and after 
filtering (Equation (\ref{eqfilter.eq})), we have, by definition of the threshold
$s(n,p,\eps)$, that $Tv=u$ with probability larger than $1-\eps$.

Now, if we want to obtain in the same way a denoising filter for white noise, we
have first to be able to relate the size threshold $s$ used for the area
openings and closings to the standard deviation $\sigma$ of the white noise
and to the grey level $\lambda$ of the thresholded image. 
In order to do this, let us consider a pure white noise image $w$:
all the $w(x)$'s are independent identically distributed random variables with
distribution $\mathcal{N}(0,\sigma^2)$
(gaussian with mean $0$ and variance $\sigma^2$).
For $\lambda\in\mathbb{R}$, let us consider the thresholded image
$w_{\lambda}=\mathbb{I}_{ w\geq \lambda}$. We then have
$$\mathbb{P}(w_{\lambda}(x)=1) = \int_{\lambda}^{+\infty}
\frac{1}{\sigma\sqrt{2\pi}} e^{-x^2/2\sigma^2} \, dx .$$
This last term, denoted by $\tilde{p}_{\lambda,\sigma}$ should be the
analogue of the probability parameter $p_\lambda$ defined in the case of
impulse noise (Equation (\ref{eqfilter.eq}) in the previous subsection). Now
the main difference 
here is that  $\tilde{p}_{\lambda,\sigma}$ is not necessarily small (it goes
to $1$ as $\lambda$ goes to $-\infty$), and the thresholded image
$w_{\lambda}$ may contain arbitrarily large connected components. Thus a
filter like the area opening or closing will never be able to remove all the
noise. The problem here is that the type of filter we have considered is not
adapted to white noise.

Generally, when using a probabilistic approach for filtering, one needs a model
for the image and one for the noise. Here, we
do not need a model for 
  the image, since we only use an ``a contrario'' hypothesis. It means that we
  only need to know 
  that ``the image is not noise'' in the sense that large connected components,
  which have a very small probability of appearing in impulse noise, necessarily belong
  to the image. This approach does not work in the case of white noise since
  the size of the connected components
of level sets is not a good way to discriminate white noise from the image
(both contain large components).
Nevertheless, if we are able to find some characteristic geometric features
(as the size of connected components in the case of impulse noise),
the proposed approach could be extended to white noise or to other models of noise.

\section{Conclusion}

We have introduced a
mathematical model for random images, in which we were able to compute the
probability of appearance of any ``local pattern'' (Theorem
\ref{th:psi}). This was then used to give an explicit formula for the size
threshold $s(n,p,\eps)$, such that the probability of appearance of a
component of size  $k\geq s(n,p,\eps)$ in a $n\times n$ image of pure
noise with probability parameter $p$ is less than $\eps$. Using this value of
$s(n,p,\eps)$ for the area openings and closings defined by Vincent will
ensure that, with probability larger than $1-\eps$, pure noise is completely
removed. This denoising process was then extended to grey level images using
their threshold decomposition.  There, the proposed
area threshold depends on both
 the probability parameter $p$ of the impulse noise and the grey level $\lambda$ of
the level set. \\
Now, some questions remain, that have not been addressed in this paper: if the
probability parameter $p$ of the impulse noise is unknown, what is the best way to
estimate it~? For a binary pure noise image, the best estimate of $p$ is
simply the ratio of the number of black pixels to the area of the image. Then,
by analogy, a first answer for binary images (like the chessboard for
example) is to compute the relative 
number of black pixels outside a dilation of the ``large'' black
components. Thus, for a grey level image, it is possible to use the threshold
decomposition to obtain initial estimates of $p_\lambda=p\times\lambda/256$
and then to estimate $p$ using, for example, a linear regression.
Now, it is not clear that this estimate will be a good one since natural grey
level images often contain textures creating small components which
over-estimate $p$. In order to obtain a reliable estimate of $p$, it would be
necessary to use also some information extracted from the statistical moments
(like the covariance, three-point probability, etc\ldots) measured on the image.

\bigskip
\bigskip

\noindent {\bf Acknowledgements}

\medskip

We would like to thank Mireille Bousquet-Mélou for all the
references about the enumeration of square lattice animals.
We also thank the anonymous referees for their detailed comments and
suggestions. 

\bibliographystyle{plain}

\bibliography{poisson}

\appendix
\section{Appendix}

\begin{demoof}{Lemma \ref{hypmoments}}

Fix $l \in \mathbb{N}^{\ast}$. Recall that $X_n$ counts the number of occurrences of the 
meaningful patterns $\bar{D}_1,\ldots,\bar{D}_{e(\psi)}$ in the random image 
$\mathcal{I}_{n,p(n)}$ where $p(n) = cn^{-\frac{2}{b(\psi)}}$. We are interested in:
$$
E_{l}(X_{n}) = \sum_{k \geq l} \mathbb{P}(X_{n} = k)\frac{k!}{(k-l)!} ~.
$$
We need to prove that $E_{l}(X_{n})$ tends to $(e(\psi)c^{b(\psi)})^{l}$ as $n$ 
tends to infinity.\\
One can see $E_{l}(X_{n})$ as the average number of ordered $l$-tuples of copies 
of the patterns $\bar{D}_{1}, \ldots , \bar{D}_{e(\psi)}$ in $\mathcal{I}_{n,p(n)}$. 
Thus, we can write:
\begin{eqnarray*}
E_{l}(X_{n}) & = & \mathbb{E} \left( \sum_{{\scriptstyle x_{1},\ldots,x_{l}}
  \atop \scriptstyle  
x_{i} \not= x_{j}} \, \, \sum_{{\scriptstyle 1 \leq j_{1},\ldots,j_{l}}  
 \leq e(\psi)} \; \mathbb{I}_{\bar{D}_{j_{1}}(x_{1}) \wedge
  \ldots \wedge \bar{D} 
_{j_{l}}(x_{l})}( \mathcal{I}_{n,p(n)}) \right) \\
& = & \sum_{s=1}^{l} \sum_{{\scriptstyle (x_{1},\ldots,x_{l})} \atop \scriptstyle 
\in {\cal C}(s)} \sum_{{\scriptstyle 1 \leq j_{1},\ldots,j_{l}} \atop \scriptstyle \leq e(\psi)} 
\mu_{n,p(n)}( \bar{D}_{j_{1}}(x_{1}) \wedge \ldots \wedge \bar{D}_{j_{l}}(x_{l}) ) 
~,
\end{eqnarray*}
where, for $s = 1,\ldots,l$, ${\cal C}(s)$ represents the set of $l-$tuples
$(x_{1},\ldots 
,x_{l})$ of pixels in $\Xi_{n}$ such that the set $\lbrace B(x_{1},r),\ldots,
B(x_{l},r)\rbrace$ is composed of $s$ equivalence classes for the $4-$connectivity relation.\\
The term corresponding to $s = l$ in the last sum will be denoted by 
$E_{l}^{'}(X_{n})$ and the rest by $E_{l}^{''}(X_{n})$. The quantity $E_{l}^{'}(X_{n})$ can 
be seen as the average number of ordered $l$-tuples of copies of $\bar{D}_{1},\ldots,\bar{D}
_{e(\psi)}$, on non-overlapping balls. We will first show that:
\begin{equation}
\label{limE'}
\lim_{n \to \infty} E_{l}^{'}(X_{n}) = (e(\psi)c^{b(\psi)})^{l} ~.
\end{equation}
Then we will prove that $E_{l}^{''}(X_{n})$ tends to $0$ as $n$ tends to infinity.
\vskip 3mm
We want to choose $l$ pixels $x_{1},\ldots,x_{l}$ such that the balls of radius $r$ 
centered on those pixels are two by two disjoint. For the first pixel $x_{1}$, there are 
$n^{2}$ possibilities. Let $2 \leq j \leq l$ and suppose pixels $x_{1},\ldots,
x_{j-1}$ have been chosen. For the $j$-th choice, the set of all pixels $x$
such that $d(x,x_{k}) \leq 2r$ for some $1\leq k\leq j-1$, must be avoided. The cardinality 
of this set is bounded by $(j-1)\times(8r^{2}+4r+1)$ whatever $x_{1},\ldots,x_{j-1}$. 
This bound does not depend on $n$. So, asymptotically the number of choices for the $j-$th 
element is $n^{2}$, and consequently the cardinality of ${\cal C}(l)$ is equivalent
to $n^{2l}$.
On the other hand, if two balls $B(x,r)$ and $B(x',r)$ are disjoint, then for all $1 
\leq j,j' \leq e(\psi)$, the random variables $\mathbb{I}_{\bar{D}_{j}(x)}$
and $\mathbb{I}_{\bar{D}_{j'}(x')}$ are independent. Therefore, we
obtain the first limit (relation (\ref{limE'})): 
$$
E_{l}^{'}(X_{n}) \sim n^{2l}(e(\psi) p(n)^{b(\psi)} (1-p(n))^{2r^{2}+2r+1-b(\psi)})^{l} \sim 
(e(\psi)c^{b(\psi)})^{l} ~.
$$
The factor $e(\psi)^{l}$ comes from the choice of the $e(\psi)$ patterns $\bar{D}_{1},\ldots,
\bar{D}_{e(\psi)}$ for the $l$ chosen balls.
\vskip 3mm
There remains to prove that $E_{l}^{''}(X_{n})$ tends to $0$ as $n$ tends to infinity. The 
intuition is that if two patterns occur in overlapping balls, then locally more than $b(\psi)$ 
black pixels are present in a ball of radius $2r$. This has vanishing probability, by Lemma 
\ref{le:threshold}.\\
Let $1 \leq s \leq l-1$ and $(x_{1},\ldots,x_{l})$ be an element of ${\cal C}(s)$. 
Let $C_{1}, \ldots ,C_{s}$ represent the
connected components  
of the set $\cup_{k=1}^l B(x_{k},r)$. Then by independence between them 
(they concern disjoint pixel sets):
$$ 
\mu_{n,p(n)}(\bar{D}_{j_{1}}(x_{1}) \wedge \ldots \wedge \bar{D}_{j_{l}}(x_{l})) = 
\prod_{m=1}^{s} \mu_{n,p(n)}( \bigwedge_{k; B(x_{k},r) \in C_{m}}
\bar{D}_{j_{k}}(x_{k}))  
~.
$$
As a consequence of $s \leq l-1$, there exists at least one connected component, say $C_{1}$, 
having at least two elements. Since the black pixel sets of two different
patterns of $\mathcal{D}_{0}(\psi)$ cannot be translated of each other, there must be at least
$b(\psi)+1$ black pixels in  
$C_{1}$. Thus we have
$$
\mu_{n,p(n)}( \bigwedge_{k; B(x_{k},r) \in C_{1}} \bar{D}_{j_{k}}(x_{k}) ) \leq p(n)^{
b(\psi)+1} ~.
$$
For the other connected components, we simply bound
$$
\mu_{n,p(n)}(\bigwedge_{k; B(x_{k},r) \in C_{m}}
\bar{D}_{j_{k}}(x_{k})) \leq \mu_{n,p(n)} 
(\bar{D}_{j_{k_{m}}}(x_{k_{m}})) \leq p(n)^{b(\psi)} ~,
$$
for any index $k_{m}$ such that $B(x_{k_{m}},r) \in C_{m}$. Therefore, we obtain the following 
result:
$$
\mu_{n,p(n)}(\bar{D}_{i_{1}}(x_{1}) \wedge \ldots \wedge
\bar{D}_{i_{l}}(x_{l}))  \leq p(n)^{sb(\psi)+1} .
$$
Finally, the set ${\cal C}(s)$ only has $O(n^{2s})$ elements 
and the number of ways to choose $l$ elements 
among $\bar{D}_{1},\ldots,\bar{D}_{e(\psi)}$ does not depend on $n$. Consequently,
the desired result follows:
$$
E_{l}^{''}(X_{n}) \leq  \sum_{s=1}^{l-1} O(n^{2s} \times
n^{-\frac{2(sb(\psi)+1)}{b(\psi)}}) = \sum_{s=1} 
^{l-1} O(n^{-\frac{2}{b(\psi)}}) = o(1) ~.
$$

\end{demoof}

\end{document}